\documentclass[a4paper,12pt]{article}
\usepackage{amsmath}
\usepackage[mathscr]{eucal}
\usepackage{amssymb}
\usepackage{latexsym}
\usepackage{amsthm}
\usepackage{amscd}
\theoremstyle{plain}
\newtheorem{theorem}{Theorem}
\newtheorem{corollary}{Corollary}[section]
\newtheorem{lemma}{Lemma}[section]
\newtheorem{proposition}{Proposition}
\theoremstyle{remark}
\newtheorem{remark}{Remark}
\theoremstyle{definition}

\begin{document}

\title{Notes on hyperelliptic fibrations of genus $3$, I\\ 
\footnotetext{During this study, the author was supported 
by DFG Forshergruppe 790
``Classification of algebraic surfaces and compact complex manifolds''.}
}
\author{Masaaki Murakami 
}
\date{}
\maketitle
\begin{abstract} 
We shall study the structure of hyperelliptic fibrations of genus $3$, 
from the view point given by Catanese and Pignatelli in \cite{fibrationsI'}.
Here by a hyperelliptic fibration of genus $3$, we mean a connected 
surjective morphism $f: S \to B$ from a nonsingular complex algebraic 
surface $S$ to a nonsingular complex projective curve $B$ with general  
fiber hyperelliptic of genus $3$.  
In this part I, we shall give a structure theorem for such fibrations 
for the case of $f : S \to B$ with all fibers $2$-connected. 
The resulting structure theorem 
is similar to one given in \cite{fibrationsI'} for genus $2$ fibrations.    
We shall also give, for the case of $B$ projective line, 
sufficient conditions for the existence of such fibrations 
$f: S \to B$'s from the view point of our structure theorem, 
prove the uniqueness of the deformation type and the simply connectedness 
of $S$ for some cases, and give some examples including 
those with simply connected $S$ and slope $3.6$ and 
those with minimal regular $S$ with geometric genus $p_g = 4$ and 
the first Chern number $c_1^2 = 8$. 
The last example turns out to be a member of the family 
$\mathcal{M}_0$ given in Bauer--Pignatelli \cite{caninvpg4c8}.   
\end{abstract}
\section{Introduction}  \label{sectn:introduction}

As is known to every algebro-geometers, the study of a fibration 
$S \to B$ of a surface $S$ over a nonsingular curve $B$ has 
been an important branch of algebraic geometry, even from  
early years of Italian school 
(\cite{complexsurf}, \cite{enriquessuperfici}). 
This is not only because a fibration itself is of an object of 
interest, but also because the study of a fibration 
of a surface gives important information on the structure of 
the surface itself.  In fact one can easily recall the 
important roles that the elliptic fibrations and ruled structures  
play in Enriquess--Kodaira classification of algebraic and analytic 
surfaces. After Kodaira's works on elliptic fibrations 
\cite{kodoncompact2}, \cite{kodoncompact3}, 
there are many works on fibrations, for which  
several approaches have been developed 
(e.g., \cite{ogggenus2}, \cite{namuengenus2}, \cite{horikawapencils}, 
\cite{xiaogenre2}, \cite{thesismlopes}, \cite{pgq1catcil}). 

One of the modern approaches for the study of fibrations 
of surfaces is that through relative canonical algebras. 
In this approach,   
the first step is the study of the canonical algebras of 
the fibers (e.g., those done for the case of the genus of the 
fibers $\leq 3$ by Mendes Lopes 
in her thesis \cite{thesismlopes}),  
since by the Krull--Azumaya lemma, the study 
of the local structures of the relative canonical algebras 
can be to some extent reduced to it. 
As one can see from for example \cite{reidpencils}, 
\cite{1-2-3}, and \cite{cliffkonno}
such study is 
important and useful for the study of global 
structures of the surfaces. 
We notice that today there are several attempts and deep results 
on the local structures of the relative canonical 
algebras and their application to the study of global structures 
of the fibrations (see \cite{akgloballocal}),  
and that as for the global structures of relative canonical algebras, 
much less are known compared to the case of the local structures,  
although we have some general theorems, 
e.g., Fujita's results on semi-positivity of the direct images of 
relative canonical sheaves
(\cite{fujitakaehler}), 
as our basic ingredients for the study in this direction. 

Recently, in \cite{fibrationsI'}, 
F.\ Catanese and R.\ Pignatelli successfully 
developed a new method for the case of genus $2$ fibrations 
and ($2$-connected) genus $3$ non-hyperelliptic firbations.     
They introduced the notion of admissible $5$-tuple, 
which is a collection of data extracted from the 
structure of the relative canonical algebra, and    
showed that the $5$-tulpe completely determines the global 
structure of the relative canonical model $X \to B$ of 
the fibration $S \to B$. 
This method turned out to be powerful. In fact, in the same 
paper, using this result, 
they were able to give a half-page proof 
of Bombieri's result on bicanonical pencils of numerical 
Godeaux surfaces, and also to prove that the moduli space 
of minimal surfaces with $c_1^2 =3$, $p_g = q = 1$, and with 
albanese fibers of genus $2$ has exactly three connected 
components (all irreducible), thus completing the classification 
of surfaces with  $c_1^2 =3$ and $p_g = q = 1$ 
initiated in \cite{pgq1catcil} and \cite{catcilsymm}
(here as usual, $c_1^2$, $p_g$, and $q$ denote the first Chern number, 
the geometric genus, and the irregularity, respectively, of a surface).

In the present paper, we study the next steps, and establish 
a theorem (similar to those in \cite{fibrationsI'}) for 
the easiest case, i.e., the case 
of hyperelliptic fibrations $S \to B$ of genus $3$ 
with all fibers $2$-connected (Theorem \ref{thm:maintheorem}). 
The resulting theorem is similar to those by Catanese 
and Pignatelli, and is the existence and description of 
one--to--one correspondence between the isomorphism classes of 
fibrations $S \to B$ as above and the isomorphism classes of 
admissible $5$-tuples, 
which will be introduced in Section \ref{scn:results}. 
We shall also give, for the case of base curve $B$ projective line, 
sufficient conditions for the existence of our fibrations 
$f: S \to B$'s from the view point of $5$-tuples
(Proposition \ref{prop:existence}), show for some cases the 
uniqueness of the deformation type and  
the simply connectedness of $S$ 
(Theorem \ref{thm:simplyconnected}, 
Proposition \ref{prop:exsistsimplconn}), 
and give some examples including 
those with topologically simply connected $S$ and slope $3.6$ 
(Remark \ref{rem:simpconnslope}) and 
those with minimal regular $S$ with $p_g =4$ and $c_1^2 = 8$
(Proposition \ref{prop:bauerpignatellim0}). 
This last example turns out to belong to the family $\mathcal{M}_0$
defined in Bauer--Pignatelli \cite{caninvpg4c8}, where they classified 
minimal regular surfaces with $p_g = 4$ and $c_1^2 = 8$ with 
canonical involution. 
The last example is of our special interest, since 
the study of surfaces with $p_g = 4$ has long history 
from Enriquess \cite{enriquessuperfici} 
(see also \cite{cilcansfpg4}), 
and after the complete classification of the case $c_1^2 = 7$
by I.\, Bauer \cite{pg4c7bauer}, the next object is the case $c_1^2 =8$ 
(see \cite{cilcansfpg4}, 
\cite{caninvpg4c8}, 
\cite{catlipign}, 
\cite{oliverioevenpg4c8}, 
\cite{supinopg4c8}).

Although the $2$-connectedness condition on the fibers is strong, 
our main result (which covers only the simplest case) 
is already useful to produce some interesting examples.  
In fact, as the first step for this purpose, 
the structure theorem for non-hyperelliptic deformations 
of the genus $3$ hyperelliptic fibrations $S \to B$ as above 
will be given in the part II of this series, 
and there 
the above mentioned surfaces with $p_g = 4$ and $c_1^2 = 8$ 
will be deformed to a family of surfaces with non-hyperelliptic genus 
$3$ fibrations. 
Here one might also notice that one of the advantages 
of treating our hyperelliptic fibrations through relative canonical 
algebras (in stead of through double cover descriptions) lies in 
that we can easily connect them to non-hyperelliptic deformations, 
as can bee seen also from the results in \cite{thesismlopes}.     
 
The present paper is organized as follows. In Section \ref{scn:results}, 
we shall  
introduce the notion of $5$-tuple and state the main theorem. 
Propositions \ref{prop:5tuplefibr} and \ref{prop:fibr5tuple} explain 
how the isomorphism classes of fibrations 
and those of $5$-tuples correspond. 
In Section \ref{scn:strrelcanalg}, 
we shall study the global structure of 
the relative canonical algebras of our fibrations $S \to B$. 
The key results are Propositions \ref{prop:amodstrofr} 
and \ref{prop:v2-descr}, 
which describe the structure of relative canonical algebras of 
our $S \to B$ (see also Remark \ref{rem:multstrdscdbydelta}).  
In Section \ref{scn:pfmaintheorem}, 
using the computation in Section \ref{scn:strrelcanalg}, 
we shall prove our main theorem, i.e., Theorem \ref{thm:maintheorem}.   
Finally, in Section \ref{scn:existenceandexs}, 
we shall give sufficient conditions for 
the existence of fibrations and study some examples. 
While Proposition \ref{prop:existence} gives sufficient conditions 
for the existence of admissible $5$-tuples, 
Theorem \ref{thm:simplyconnected},  
Propositions \ref{prop:exsistsimplconn} and \ref{prop:bauerpignatellim0}, 
and Remarks \ref{rem:someexamples} and \ref{rem:simpconnslope}, 
etc., study examples. 
   
Throughout this paper, we work over the complex number field 
$\mathbb{C}$. 

\medskip
{\sc Acknowledgment}
\medskip

\noindent
The author expresses deepest gratitude to Prof.\ Fabrizio Catanese 
and Prof.\ Ingrid Bauer 
for the comfortable environment for the study, financial support, 
discussions, and also for welcoming him for the study in this direction. 
During this study, the author was supported by DFG Forshergruppe 790
``Classification of algebraic surfaces and compact complex manifolds''.  

\medskip
{\sc Notation}
\medskip

\noindent
In this article, the symbol $k$ always denotes the complex 
number field $\mathbb{C}$.  
If $V$ is a locally free sheaf of finite rank on a scheme, 
$\mathrm{rk} \, V$ denotes its rank. 
The symbols $\mathcal{S} (V)$ and $\mathrm{Sym}^n V$ denote 
the symmetric tensor algebra (of $V$) and its homogeneous part 
of degree $n$, respectively. 
If $R$ is a graded algebra, $R_j$ denotes its homogeneous part 
of degree $j$. Thus, for example, for the polynomial ring 
$k[x_0, x_1, x_2]$ over the complex number field $k= \mathbb{C}$, 
an element in $k[x_0, x_1, x_2]_j$ is 
a homogeneous polynomial in $x_0$, $x_1$, $x_2$ of degree $j$. 
The symbol $\amalg$ means taking the disjoint union of sets. 
If $p$ is a point of a scheme, $k(p)$ denotes the residue field 
at $p$ of this scheme. If $\mathcal{F}$ is a sheaf on a scheme, 
$\mathcal{F}_p$ denotes the stalk at $p$ of $\mathcal{F}$. 
If moreover the scheme is over $k$, and $\mathcal{F}$ is 
coherent, $h^i (\mathcal{F})$ denotes the dimension over $k$ 
of the $i$-th cohomology group $H^i (\mathcal{F})$ of $\mathcal{F}$.  
If $S$ is a scheme, and $D$, a Cartier divisor on $S$, then 
$\mathcal{O}_S$ and $\mathcal{O}_S (D)$ denote the structure sheaf 
of $S$ and the invertible sheaf associated to $D$, respectively. 
If $S$ is projective and smooth over $k$, the symbol $K_S$ 
as usual denotes the canonical divisor of $S$. 

\section{Statement of the main result}  \label{scn:results}

In this section, we shall state the main theorem. 
The result is the existence of one--to--one correspondence 
between isomorphism classes of a kind of data, 
which we shall call $5$-tuples, and isomorphism classes of 
genus $3$ hyperelliptic fibrations with all fibers $2$-connected. 
In order to state the result, we first introduce the notion of 
$5$-tuple, and observe how to associate a fibration to 
a $5$-tuple, and then a $5$-tuple to a fibration. 
And then we state the main theorem.

\subsection{From $5$-tuples to fibrations}  \label{subscn:5tofibr}

First, let us define the $5$-tuple and observe how to 
associate a genus $3$ fibration to it. 
Let $B$ be a smooth projective curve over the complex number field 
$k = \mathbb{C}$. 
Let $V_1$ and $V_2^+$ be locally free sheaves on $B$ of rank 
$\mathrm{rk}\, V_1 = 3$ and $\mathrm{rk}\, V_2^+ = 5$, respectively.  
 
Assume that we are given a surjective morphism 
$\sigma_2 : \mathrm{Sym}^2 V_1 \to V_2^+$ of sheaves. 
Then the kernel $L = \ker \sigma_2$ of $\sigma_2$ is invertible, 
and for each natural number $n \in \mathbb{N}$ the natural inclusion 
$L \to \mathrm{Sym}^2 V_1$ induces an injective morphism 
$L \otimes \mathrm{Sym}^{n-2} V_1 \to \mathrm{Sym}^n V_1$. 
So we define the coherent $\mathcal{O}_B$--module $\mathcal{A}_n$
by the following short exact sequence:
\[
0 \to L \otimes \mathrm{Sym}^{n-2} V_1 \to \mathrm{Sym}^n V_1
\to \mathcal{A}_n \to 0. 
\]
Then $\mathcal{A}_n$ is a locally free sheaf of 
rank $2n + 1$. 
Let $\mathcal{S} (V_1)$ be the symmetric $\mathcal{O}_B$--algebra
associated to the  $\mathcal{O}_B$--module $V_1$, and 
put $\mathcal{A} = \bigoplus_{n=0}^{\infty} \mathcal{A}_n$
Then, via the natural projection 
$\mathcal{S} (V_1) \to \mathcal{A}$, 
the algebra structure of 
$\mathcal{S} (V_1)$ induces a 
quasi-coherent graded 
$\mathcal{O}_B$--algebra structure 
on the direct sum $\mathcal{A} = \bigoplus_{n=0}^{\infty} \mathcal{A}_n$.  

Using these graded algebras, 
we define the two varieties $\mathcal{C}$ and $\mathbb{P}$ by 
$\mathcal{C} = \mathcal{P}\mathrm{roj}\, \mathcal{A}$ and 
$\mathbb{P} = \mathcal{P}\mathrm{roj}\, \mathcal{S}(V_1) = \mathbb{P} (V_1)$. 
By the projection $\mathcal{S}\, (V_1) \to \mathcal{A}$, 
we obtain a natural closed embedding $\mathcal{C} \to \mathbb{P}$
over the curve $B$.  

Set $V_2^- = (\det V_1 ) \otimes L^{\otimes (-1)}$ and 
$\mathcal{R} = \mathcal{A} \oplus (\mathcal{A}[-2] \otimes V_2^-)$, 
where $\mathcal{A} [-2]$ is the $(-2)$--shift of the graded algebra
$\mathcal{A}$. 
Then the $\mathcal{O}_B$--module $\mathcal{R}$ allows a natural 
graded $\mathcal{A}$-module structure.
Assume moreover that we are given an element
\begin{align}
\delta \in \mathrm{Hom}_{\mathcal{O}_B}
((V_2^-)^{\otimes 2}, \, \mathcal{A}_4)
&\simeq 
H^0 (B, \, \mathcal{A}_4 \otimes  (V_2^-)^{\otimes (-2)} ) \notag \\
&\simeq 
H^0 (\mathcal{C}, \mathcal{O}_{\mathcal{C}} (4) 
\otimes \pi_{\mathcal{C}}^* {(V_2^-)^{\otimes (-2)}}), \notag
\end{align}
where $\pi_{\mathcal{C}} : \mathcal{C} \to B$ is the 
natural projection. 

Since $\delta : (V_2^-)^{\otimes 2} \to \mathcal{A}_4$ 
induces a natural morphism of $\mathcal{O}_B$--modules 
$ (\mathcal{A}[-2] \otimes V_2^-)_m \otimes 
(\mathcal{A}[-2] \otimes V_2^-)_n   
\to \mathcal{A}_{m+n}$ for each $m$, $n \geq 2$, 
and since the graded $\mathcal{A}$--module structure on $\mathcal{R}$ 
gives $\mathcal{R}_m \otimes \mathcal{A}_n 
\to \mathcal{R}_{m+n}$, 
the element $\delta$ determines a graded $\mathcal{O}_B$--algebra
structure on the $\mathcal{A}$--module $\mathcal{R} = 
\mathcal{A} \oplus (\mathcal{A}[-2] \otimes V_2^-)$.
Note that the natural inclusion $\mathcal{A} \to 
\mathcal{R} = \mathcal{A} \oplus (\mathcal{A}[-2] \otimes V_2^-)$ 
of the first direct summand $\mathcal{A}$ 
is a morphism of graded $\mathcal{O}_B$--algebras. 
Thus if we put $X = \mathcal{P}\mathrm{roj}\, \mathcal{R}$, 
we obtain the following commutative diagram: 

\begin{equation*} 
\begin{CD}  
   X =  \mathcal{P}\mathrm{roj}\, \mathcal{R}  @>\text{$$}>> \mathcal{C} = \mathcal{P}\mathrm{roj}\, \mathcal{A} @>\text{$$}>> 
   \mathbb{P} = \mathcal{P}\mathrm{roj}\, \mathcal{S} (V_1)\\
   @VV\text{$\bar{f}$}V  @VV\text{$\pi_{\mathcal{C}}$}V @VV\text{$\pi_{\mathbb{P}}$}V\\
   B    @>\text{$=$}>> B @>\text{$=$}>> B \, , 
\end{CD}
\end{equation*}
where $\bar{f}$, $\pi_{\mathcal{C}}$, and $\pi_{\mathbb{P}}$ are 
the natural projections. 

We shall call 
$(B,  V_1, V_2^+, \sigma_2, \delta)$ as above 
a $5$-tuple for relatively minimal ($2$-connected) hyperelliptic 
fibrations of genus $3$, or simply, a $5$-tuple.  
We shall call $\bar{f}: X \to B$ as above the relative canonical model
associated to the $5$-tuple, or simply, the associated relative canonical 
model. If the associated relative canonical model 
$\bar{f} : X \to B$ satisfies the 
two conditions 
\smallskip

I) $\mathcal{C}$ has at most rational double points as its singularities, and 

II) $X$ has at most rational double points as its singularities, 
\smallskip

\noindent  
we say that the $5$-tuple 
$(B,  V_1, V_2^+, \sigma_2, \delta)$ is admissible. 

We shall prove the following proposition in Section \ref{scn:pfmaintheorem}.

\begin{proposition}    \label{prop:5tuplefibr}
Let $(B,  V_1, V_2^+, \sigma_2, \delta)$ be 
an admissible $5$-tuple, 
$\bar{f} : X \to B$, its associated relative canonical model,
and $S \to X$, the minimal resolution of singularities of $X$. 
Denote by $f: S \to B$ the composite of the two morphisms 
$S \to X$ and $\bar{f} : X \to B$. 
Then $f: S \to B$ is a relatively minimal hyperelliptic fibration 
of genus $3$ all of whose fibers are $2$-connected.  
\end{proposition}

\subsection{From fibrations to $5$-tuples}

Next, let us observe how to associate a $5$-tuple to 
a hyperelliptic fibration of genus $3$. 
Let $f: S \to B$ be a relatively minimal hyperelliptic fibration  
of genus $3$ 
all of whose fibers are $2$-connected. 
Let $\omega_{S | B}$ be the relative dualizing sheaf of $f$,  
and $V_n = f_* (\omega_{S | B}^{\otimes n})$, the direct image sheaf 
of $\omega_{S | B}^{\otimes n}$ by $f$.   
Then we have a natural decomposition 
$V_n = V_n^+ \oplus V_n^-$ into eigen-sheaves: 
$V_n^+$ and $V_n^-$ are eigen-sheaves of eigenvalue $+1$ and 
$-1$, respectively, with respect to the action by the hyperelliptic 
involution of $f$. 
It is easy to see that $\mathrm{rk}\, V_1 (=V_1^-) =3$ and 
$\mathrm{rk}\, V_2^+ = 5$.
Now let $\sigma_2 : \mathrm{Sym}^2 V_1 \to V_2^+$ be the natural 
morphism induced by the multiplication structure of the 
relative canonical algebra 
$\mathcal{R} = \bigoplus_{n=0}^{\infty} V_n$ of $f$, 
and $L = \mathrm{ker}\, \sigma_2$, its kernel.  
We denote by $\delta : (V_2^-)^{\otimes 2} \to V_4^+$ the 
natural morphism induced by the multiplication structure of $\mathcal{R}$. 

We shall prove the following proposition in Section \ref{scn:pfmaintheorem}.

\begin{proposition} \label{prop:fibr5tuple}
Let $f: S \to B$ a relatively minimal hyperelliptic fibration 
of genus $3$ 
all of whose fibers are $2$-connected,  
$\mathcal{R} = \bigoplus_{n=0}^{\infty} V_n$, 
its relative canonical algebra,  
and $\mathcal{A} \subset \mathcal{R}$, 
the graded $\mathcal{O}_B$--subalgebra 
generated by the degree $1$ part $V_1 = V_1^-$ in $\mathcal{R}$.  
Then the degree $4$ part $\mathcal{A}_4$ of $\mathcal{A}$ coincides 
with $V_4^+$, and $(B, V_1, V_2^+, \sigma_2, \delta)$ above 
forms an admissible $5$-tuple for relatively minimal 
hyperelliptic fibrations of genus $3$. 
\end{proposition} 

Given a fibration $f: S \to B$ as above, 
we shall call $(B, V_1, V_2^+, \sigma_2, \delta)$ 
as in Proposition \ref{prop:fibr5tuple} the $5$-tuple associated to $f$. 

\subsection{Main theorem}

Under the terminology as above, our main theorem, which we 
shall prove in Section \ref{scn:pfmaintheorem}, is the following: 

\begin{theorem}  \label{thm:maintheorem}
Let $B$ be a smooth projective curve over a complex number field 
$\mathbb{C}$. 
Then via the associations given in 
Propositions \ref{prop:5tuplefibr} and \ref{prop:fibr5tuple},
which are mutually inverse,  
the isomorphism classes of relatively minimal genus $3$ hyperelliptic 
fibrations  
with all fibers $2$-connected are in one-to-one correspondence 
with the isomorphism classes of admissible $5$-tuples 
$(B, V_1, V_2^+, \sigma_2, \delta)$'s.   

Moreover, given an admissible $5$-tuple 
$(B, V_1, V_2^+, \sigma_2, \delta)$, the resulting 
surface $S$ appearing in the associated fibration $f: S \to B$ 
has numerical invariants

\begin{align}
\chi (\mathcal{O}_S) &= \deg V_1 + 2 (b-1) \notag \\
c_1^2 (S) & = 4 \deg V_1 - 2 \deg L + 16 (b-1)  \notag 
\end{align} 
where $L$ is the kernel of the morphism $\sigma_2$, 
and $b = g(B)$, the genus of the base curve $B$. 

\end{theorem}

\section{The structure of relative canonical algebra}  
\label{scn:strrelcanalg}

In this section, we shall study the structure of the relative 
canonical algebras for our fibrations. 
Let $f: S \to B$ be a relatively minimal hyperelliptic genus $3$
fibration with all fibers $2$-connected. 
Recall that we have a natural decomposition 
$V_n = V_n^+ \oplus V_n^-$ of the direct image sheaf 
$V_n = f_* (\omega_{S|B}^{\otimes n})$ induced by the action 
 on $V_n$ by the hyperelliptic involution of $f$. 
It is straightforward to see that both $V_n^+$ and $V_n^-$ are 
locally free sheaves and to see that we have  
\begin{align}
\mathrm{rk} V_n^{\pm} &= 2n + 1 \notag \\
\mathrm{rk} V_n^{\mp} &= 2n - 3, \notag
\end{align}  
where the symbol $\pm$ stands for $ +$ if $n$ is even, for $-$ if $n$ is odd, 
and the symbol $\mp$ stands for $-$ if $n$ is even, for $+$ if $n$ is odd. 
In using the symbols $\pm$ and $\mp$, we shall keep this rule throughout 
this paper. 

Note that we have in particular
\[
 \mathrm{rk}\, V_1 = \mathrm{rk}\, V_1^{-} = 3, \qquad
\mathrm{rk}\, V_2^{+} = 5, \qquad
\mathrm{rk}\, V_2^{-} = 1. 
\]

\begin{lemma}   \label{lm:canringfiber}
Let $F$ be a fiber of the fibration $f$, and 
$R(F, K_F)$, the canonical ring of the fiber $F$. 
Then 
\[
R(F, K_F) \simeq k[x_0, x_1, x_2, y] / (Q,\, y^2 - P)
\]
as graded $k=\mathbb{C}$-algebras, where 
$\deg x_i = 1$ for $0 \leq i \leq 2$, 
$\deg y =2$, and 
\begin{align}
Q &= x_2^2 - Q_1 (x_0, x_1) x_2 - Q_2 (x_0, x_1)
\in k[x_0, x_1, x_2]_2 \notag \\
P &= P_3 (x_0, x_1) x_2 + P_4(x_0, x_1) 
\in k[x_0, x_1, x_2]_4.  \notag
\end{align} 
Here 
$Q_j(x_0, x_1) \in k[x_0, x_1]_j$ $($ $j= 1, 2$ $)$ and
$P_i(x_0, x_1) \in k[x_0, x_1]_i$ $($ $i= 3, 4$ $)$ are 
homogeneous polynomials in $x_0$, $x_1$ of 
degree $j$ and $i$, respectively. 
\end{lemma}

Proof. See Mendes Lopes  
\cite[Theorem $6$.$1$, p.198]{thesismlopes}. 
Note that in our case 
the fiber $F$ is hyperelliptic. \qed

\begin{remark}  \label{rem:yisabase}
Note that in Lemma \ref{lm:canringfiber} the element $y$ is a 
base of $1$-dimensional linear space $R(F, K_F)_2^-$,  
i.e., that of the eigenspace of the eigenvalue $-1$ 
with respect to the action by the hyperelliptic 
involution on the homogeneous part $R(F, K_F)_2$ of degree $2$. 
The action by the hyperelliptic involution on $R(F, K_F)$ 
is given by 
\[
(x_0, x_1, x_2, y) \mapsto (-x_0, -x_1, -x_2, -y). 
\]
\end{remark}

From Lemma \ref{lm:canringfiber} we infer the 
following:

\begin{lemma}   \label{lm:baseofrfkf}
Let $(*)_n$ and $(**)_n$ be the two sets of monomials 
of weighted degree $n$ 
in $x_0$, $x_1$, $x_2$, and $y$ 
defined as follows$:$ 
\begin{align}
(*)_n = &\{ (\textrm{ monomials in  $x_0$, $x_1$ of degree $n$}) \} \notag \\
        &\amalg 
         \{ (\textrm{ monomials in $x_0$, $x_1$ of degree $n-1$}) 
         \times x_2 \}  \notag \\
(**)_n = &\{ (\textrm{ monomials in  $x_0$, $x_1$ of degree $n-2$}) 
         \times y \} \notag \\
        &\amalg 
         \{ (\textrm{ monomials in $x_0$, $x_1$ of degree $n-3$}) 
         \times x_2 y \}. \notag
\end{align}
Then for any integer $n \geq 2$, the set $(*)_n$ 
forms a base of $R(F, K_F)_n^{\pm}$, and 
the set $(**)_n$, a base of $R(F, K_F)_n^{\mp}$, where $R(F, K_F)_n^+$
and $R(F, K_F)_n^-$ are eigenspaces of eigenvalues $1$ and 
$-1$, respectively, with respect to the action by the hyperelliptic 
involution on the homogeneous part $R(F, K_F)_n$ of degree $n$.  
Here as before 
the symbol $\pm$ stands for $+$ if $n$ is even, for $-$ if $n$ is odd, 
and 
the symbol $\mp$ stands for $-$ if $n$ is even, for $+$ if $n$ is odd.  
\end{lemma}

Proof. Note that $y^2 - P$ is monic in $y$ and that $Q$ is monic 
in $x_2$. Thus the assertion follows immediately 
from Lemma \ref{lm:canringfiber}. 
\qed 

Now let us denote by $\sigma_n$ the natural morphism 
$\sigma_n : \mathrm{Sym}^n V_1 \to V_n^{\pm}$ determined by 
the multiplication structure of the relative canonical algebra 
$\mathcal{R} (f) = \bigoplus_{n = 0}^{\infty} V_n$. 
Let us moreover denote by $L= \mathrm{ker}\, \sigma_2$ the 
kernel of the morphism $\sigma_2: \mathrm{Sym}^2 V_1 \to V_2^+$. 
Then since $\sigma_2 \otimes k(p) : 
(\mathrm{Sym}^2 V_1 ) \otimes k(p) \to V_2^+ \otimes k(p)$  
is surjective for any point $p = f(F)$ by Lemma \ref{lm:baseofrfkf}, 
we obtain by the Krull--Azumaya Lemma (i.e., Nakayama's Lemma) 
the exactness of the complex 
\[
0 \to L \to \mathrm{Sym}^2 V_1 \to V_2^+ \to 0 .
\] 
Since we have $\mathrm{rk}\, \mathrm{Sym}^2 V_1 = 6$ and 
$\mathrm{rk}\, V_2^+ = 5$, we see that the sheaf $L$ is 
invertible on $B$.  

The next two lemmas follow immediately from the Krull--Azumaya Lemma:

%
%
%
%

%
%
\begin{lemma}   \label{lm:canalginnbd}
Let $F$ be a (closed) fiber of $f$, and $P$ and $Q$, 
polynomials as in Lemma \ref{lm:canringfiber} 
in $x_0$, $x_1$, and $x_2$.
Then there exists a neighbourhood $U$ of the point $p = f(F) \in B$, 
such that the relations $Q$ and $y^2 - P$ in $R(F, K_F)$ lift to 
the relations $\tilde{Q}$ and $y^2 - \tilde{P}$, respectively,   
of the form
\begin{align}
\tilde{Q} &= x_2^2 - \tilde{Q}_1 (x_0, x_1) x_2 - \tilde{Q}_2 (x_0, x_1)
    \notag \\
y^2 - \tilde{P} &= y^2 - \tilde{P}_3 (x_0, x_1) x_2 - \tilde{P}_4 (x_0, x_1)
    \notag
\end{align}
in the relative canonical 
algebra $\mathcal{R} = \mathcal{R} (f)$. 
Here $\tilde{Q}_j (x_0, x_1) \in \mathcal{O}_B (U) [x_0, x_1]_j$ 
$($$j =1, 2$$)$ and 
$\tilde{P}_i (x_0, x_1) \in \mathcal{O}_B (U) [x_0, x_1]_i$ 
$($$i=3, 4$$)$ are homogeneous polynomials 
with coefficients in $\mathcal{O}_B (U)$
of degree $j$ and $i$, respectively,  
and satisfy $\tilde{Q}_j |_{f(F)} = Q_j$ and $\tilde{P}_i |_{f(F)} = P_i$.  
\end{lemma}

\begin{lemma}  \label{lm:surjsigman}
Let $f: S \to B$ be a fibration as in the beginning of this section. 
Then for any integer $n \geq 2$, the morphism 
$\sigma_n : \mathrm{Sym}^n V_1 \to V_n^{\pm}$ is surjective. 
\end{lemma}

In fact, take liftings to the stalk $V_{1,\, p}$ 
at $p = f(F)$ of the elements
$x_0$, $x_1$, and $x_2 \in V_{1,\, p} \otimes k(p)$ 
in Lemma \ref{lm:canringfiber}. 
We use the same symbols $x_0$, $x_1$, and $x_2$ for the 
respective liftings to the stalk $V_{1,\, p}$.  
Then the polynomial $Q$ in Lemma \ref{lm:canringfiber} defines 
an element of the stalk $V_{2, \, p}^+$, 
for which we use the same symbol $Q$.   
By the Krull--Azumaya Lemma, the stalk $V_{2,\, p}^+$ at $p$ 
of the sheaf $V_2^+$ 
is generated by the elements in $(*)_2$, and since 
$Q \mapsto 0$ by the natural projection  
$V_{2,\, p}^+ \to V_{2,\, p}^+ \otimes k(p)$, we see that 
the polynomial $Q$ is in $V_{2,\, p}^+$ a linear combination 
of elements in the set $(*)_2$ with coefficients in 
$\mathfrak{M}_p$, 
where $\mathfrak{M}_p$ is the maximal ideal 
of the local ring $\mathcal{O}_{B,\, p}$ at $p$. 
This shows the existence of $\tilde{Q}$
in Lemma \ref{lm:canalginnbd}. 
The proof for the existence 
of $\tilde{P}$ in Lemma \ref{lm:canalginnbd} is the same. 
Lemma \ref{lm:surjsigman} can be proved by the same 
argument as that for the surjectivity of $\sigma_2$.  

Now let us define the graded $\mathcal{O}_B$--algebra 
$\mathcal{A} = \bigoplus_{n=0}^{\infty} \mathcal{A}_n$ as 
the $\mathcal{O}_B$--subalgebra 
generated by $V_1 = V_1^-$ in the relative canonical algebra 
$\mathcal{R} = \mathcal{R} (f) = \bigoplus_{n=0}^{\infty} V_n$.
By Lemma \ref{lm:surjsigman}, we see that 
\[
\mathcal{A}_n = V_n^{\pm} . 
\]

The next lemma says that the graded $\mathcal{O}_B$--algebra 
structure of $\mathcal{A}$ is completely determined by 
the natural inclusion $L \to \mathrm{Sym}^2 V_1$ 
(and hence by the projection $\sigma_2 : \mathrm{Sym}^2 V_1 \to V_2^+$). 
\begin{lemma} \label{lm:dscran}
Consider the natural morphism 
\[
L \otimes \mathrm{Sym}^{n-2} V_1  \to \mathrm{Sym}^n V_1   
:  l \otimes q \mapsto \, \tilde{Q} q 
\]
induced by the inclusion $L \to \mathrm{Sym}^2 V_1: l \mapsto \tilde{Q}$, 
where $q$ is a local section to $\mathrm{Sym}^{n-2} V_1$, and 
$l$, the local base of $L$ corresponding to the local section  
$\tilde{Q}$ to $\mathrm{Sym}^2 V_1$ 
as in Lemma \ref{lm:canalginnbd}. Then the natural complex 
\begin{equation} \label{eqn:exactseqdscran}
0 \to  L \otimes \mathrm{Sym}^{n-2} V_1\to \mathrm{Sym}^n V_1
  \to \mathcal{A}_n \to 0
\end{equation}
is exact. Moreover for any closed point 
$p \in B$, the complex (\ref{eqn:exactseqdscran}) tensored 
by the residue field $k(p)$ is also exact. 
\end{lemma}

Proof. 
Since we have 
$\tilde{Q}|_p = Q = x_2^2 - x_2 Q_1(x_0, x_1) - Q_2 (x_0, x_1)$, 
the morphism 
$(L \otimes \mathrm{Sym}^{n-2} V_1) \otimes k(p)
\to (\mathrm{Sym}^n V_1) \otimes k(p)$ is injective and has 
maximal rank. But we have
\begin{align}
&\mathrm{rk}\, (\mathrm{Sym}^n V_1) = (n+2)(n+1)/2, \notag \\
&\mathrm{rk}\, (L \otimes \mathrm{Sym}^{n-2} V_1) = n (n-1) / 2, \notag \\
&\mathrm{rk}\, \mathcal{A}_n = \mathrm{rk}\, Vn^{\pm} = 2n +1, \notag 
\end{align}
and hence 
$\mathrm{rk}\, (\mathrm{Sym}^n V_1) - \mathrm{rk}\, (L \otimes \mathrm{Sym}^{n-2} V_1) = \mathrm{rk}\, \mathcal{A}_n$.
Thus the complex (\ref{eqn:exactseqdscran}) 
tensored by $k(p)$ is exact at every closed point 
$p \in B$. This implies in particular the surjectivity of 
$L \otimes \mathrm{Sym}^{n-2} V_1\to \mathrm{ker}\, \sigma_n$. 
Since two locally free sheaves $L \otimes \mathrm{Sym}^{n-2} V_1$ and 
$\mathrm{ker}\, \sigma_n$ have the same rank $n (n-1) /2$, this shows 
that $L \otimes \mathrm{Sym}^{n-2} V_1\to \mathrm{ker}\, \sigma_n$ 
is an isomorphism, hence the assertion. \qed 
 
\begin{lemma}    \label{lm:vnmp}
The morphism 
\begin{equation}   \label{eqn:anm2v2misomvnmp}
\mathcal{A}_{n-2} \otimes V_2^- \to V_n^{\mp}: \quad q \otimes y \mapsto q y
\end{equation}
is an isomorphism.
\end{lemma}

Proof. By Lemma \ref{lm:baseofrfkf}, the morphism 
$(\mathcal{A}_{n-2} \otimes V_2^-)  \otimes k(p)  \to 
V_n^{\mp} \otimes k(p)$ is an isomorphism for any closed $p \in B$. 
Thus by the Krull--Azumaya Lemma, 
$\mathcal{A}_{n-2} \otimes V_2^- \to V_n^{\mp}$ is surjective. 
But both hands of this morphism are 
locally free sheaves of the same rank $2n-3$, 
hence the assertion. \qed

Thus we obtain the following proposition:

\begin{proposition} \label{prop:amodstrofr}
Let $\mathcal{A} \subset \mathcal{R}$ be the graded 
subalgebra as above. Then 
\[
\mathcal{R} \simeq \mathcal{A} \oplus (\mathcal{A} [-2] \otimes V_2^- )
\] 
as graded $\mathcal{A}$--modules. 
\end{proposition}

We denote by 
$\delta \in \mathrm{Hom}_{\mathcal{O}_B} ((V_2^-)^{\otimes 2}, \, \mathcal{A}_4)$ 
the $\mathcal{O}_B$--module homomorphism 
from 
$(V_2^-)^{\otimes 2}$  
to 
$\mathcal{A}_4 = V_4^+$ determined by the multiplication structure 
of the relative canonical algebra $\mathcal{R}$ 
(notice that $V_2^- \subset \mathcal{R}_2$ and 
$\mathcal{A}_4 = V_4^+ \subset \mathcal{R}_4$).  

\begin{remark}   \label{rem:multstrdscdbydelta}
By Remark \ref{rem:yisabase} and 
Lemma \ref{lm:canalginnbd}, the section $y$ is a local base 
of the invertible sheaf $V_2^-$. 
Note that with this $y$, 
the morphism $\delta : (V_2^-)^{\otimes 2} \to \mathcal{A}_4$
is given by 
$y^2 \mapsto \tilde{P} = \tilde{P}_3 (x_0, x_1) x_2 + \tilde{P}_4 (x_0, x_1)$ 
where $\tilde{P}$, $\tilde{P}_3$, and $\tilde{P}_4$ are polynomials as in 
Lemma \ref{lm:canalginnbd}. 
Moreover, via the isomorphism (\ref{eqn:anm2v2misomvnmp}) 
in Lemma \ref{lm:vnmp}, the 
multiplication morphism 
$V_m^{\mp} \otimes V_n^{\mp} \to V_{m+n}^{\pm} = \mathcal{A}_{m+n}$ 
is given by 
\[
(\mathcal{A}_{m-2} \otimes V_2^-) \otimes (\mathcal{A}_{n-2} \otimes V_2^-)
\to \mathcal{A}_{m+n} : \,  
(\alpha \otimes y) \otimes (\beta \otimes y) \mapsto 
  (\alpha \beta) \delta (y^2). 
\]
From this together with the short 
exact sequence (\ref{eqn:exactseqdscran}) 
in Lemma \ref{lm:dscran}, 
we see that the $\mathcal{O}_B$--algebra structure 
of the relative canonical algebra 
$\mathcal{R} \simeq \mathcal{A} \oplus (\mathcal{A} [-2] \otimes V_2^- )$ 
is completely described by the two morphisms $\sigma_2$ and $\delta$.   
We shall use this fact later in the proof of our main theorem. 
\end{remark}

For the $\mathcal{O}_B$--algebras above, let us denote 
by $\pi_{\mathcal{C}}: \mathcal{C} 
= \mathcal{P}\mathrm{roj}\, \mathcal{A} \to B$ and 
$\bar{f}: X 
= \mathcal{P}\mathrm{roj}\, \mathcal{R} \to B$ the 
structure morphisms of $\mathcal{C} = \mathcal{P}\mathrm{roj}\, \mathcal{A}$ 
and $X = \mathcal{P}\mathrm{roj}\, \mathcal{R}$, respectively. 
Then by the two natural morphisms of $\mathcal{O}_B$--algebras 
$\mathcal{S} (V_1) \to \mathcal{A}$ and 
$\mathcal{A} \to \mathcal{R} \simeq 
\mathcal{A} \oplus (\mathcal{A} [-2] \otimes V_2^-)$, 
we obtain the commutative diagram

\begin{equation*} 
\begin{CD}  
   X =  \mathcal{P}\mathrm{roj}\, \mathcal{R}  @>\text{$\psi$}>> \mathcal{C} = \mathcal{P}\mathrm{roj}\, \mathcal{A} @>\text{$$}>> 
   \mathbb{P} = \mathcal{P}\mathrm{roj}\, \mathcal{S} (V_1)\\
   @VV\text{$\bar{f}$}V  @VV\text{$\pi_{\mathcal{C}}$}V @VV\text{$\pi_{\mathbb{P}}$}V\\
   B    @>\text{$=$}>> B @>\text{$=$}>> B \, , 
\end{CD}
\end{equation*}
where $\psi$ is the natural projection induced by 
$\mathcal{A} \to \mathcal{R}$ above.  

Next let us describe $V_2^-$ using the locally free sheaves $V_1$ and $L$. 

\begin{lemma} \label{lm:classofc}
$\mathcal{C}$ is a divisor on $\mathbb{P}$, 
and is a member of the linear system 
$| \mathcal{O}_{\mathbb{P}} (2) \otimes \pi_{\mathbb{P}}^* (L^{\otimes (-1)})|$.
\end{lemma}

Proof. 
By (\ref{eqn:exactseqdscran}) in 
Lemma \ref{lm:dscran}, we obtain the short exact sequence
\[
0 \to \mathcal{O}_{\mathbb{P}} (-2) \otimes \pi_{\mathbb{P}}^* L 
  \to \mathcal{O}_{\mathbb{P}}
  \to \mathcal{O}_{\mathcal{C}}  \to 0 . 
\]
From this we infer  
$ \mathcal{C} \in 
| \mathcal{O}_{\mathbb{P}} (2) \otimes \pi_{\mathbb{P}}^* (L^{\otimes (-1)})|$, 
hence the assertion. \qed

\begin{lemma} \label{lm:omegacb}
Put $\mathcal{M} = (\mathrm{det}\, V_1 ) \otimes L^{\otimes (-1)}$. 
Then the following hold$:$

i$)$ $\omega_{\mathcal{C} | B} \simeq 
   \mathcal{O}_{\mathcal{C}} (-1) \otimes \pi_{\mathcal{C}}^* (\mathcal{M})$$;$ 

ii$)$ $\psi_* \mathcal{O}_X \simeq 
\mathcal{O}_C \oplus 
(\mathcal{O}_C (-2) \otimes \pi_{\mathcal{C}}^* (V_2^-))$$;$ 

iii$)$ $\psi_* (\omega_{X|B}) \simeq 
(\mathcal{O}_{\mathcal{C}}(-1) \otimes \pi_{\mathcal{C}}^* (\mathcal{M}))
\oplus
(\mathcal{O}_{\mathcal{C}}(1) \otimes 
\pi_{\mathcal{C}}^* (\mathcal{M} \otimes (V_2^{-})^{\otimes (-1)}))
$. 
\end{lemma}

Proof. 
The assertion i) follows from the adjunction formula together with 
Lemma \ref{lm:classofc} and 
$\omega_{\mathbb{P}} 
\simeq \mathcal{O}_{\mathbb{P}} (-3) \otimes 
\pi_{\mathbb{P}}^* (\omega_B \otimes \mathrm{det}\, V_1)$. 
The assertion ii) follows from the fact that 
$\psi : X = \mathcal{P}\mathrm{roj} \mathcal{R} \to 
\mathcal{C} = \mathcal{P}\mathrm{roj} \mathcal{A}$ is 
the morphism induced by  
$\mathcal{A} \to \mathcal{R} \simeq 
\mathcal{A} \oplus (\mathcal{A} [-2] \otimes V_2^-)$. 
Thus we only need to show the assertion iii). 
But since $\psi : X \to \mathcal{C}$ is finite, we have 
$\omega_X \simeq f^{!} \omega_{\mathcal{C}}$, i.e., 
$\psi_* \omega_X \simeq 
\mathcal{H}\mathrm{om}_{\mathcal{O}_{\mathcal{C}}} 
(\psi_* \mathcal{O}_X ,\, \omega_{\mathcal{C}})$, 
and hence 
$\psi_* (\omega_{X | B}) \simeq 
\mathcal{H}\mathrm{om}_{\mathcal{O}_{\mathcal{C}}} 
(\psi_* \mathcal{O}_X ,\, \omega_{\mathcal{C} | B})$. 
Then the assertion iii) follows from i) and ii). \qed 

\begin{lemma}     \label{lm:omegaxb}
There exists a natural isomorphism 
\[
\omega_{X | B} \to 
\psi^* 
(
\mathcal{O}_{\mathcal{C}}(1) \otimes
\pi_{\mathcal{C}}^* (\mathcal{M} \otimes (V_2^{-})^{\otimes (-1)})
).
\]
The inverse of this isomorphism is induced by 
the natural inclusion  
\begin{equation} \label{eqn:incltodrctim}
\mathcal{O}_{\mathcal{C}}(1) \otimes 
\pi_{\mathcal{C}}^* (\mathcal{M} \otimes (V_2^{-})^{\otimes (-1)})
\to
\psi_* (\omega_{X|B})  
\end{equation}
of the second direct summand of iii$)$, Lemma \ref{lm:omegacb}.
\end{lemma}

Proof. 
As we have shown in the proof of Lemma \ref{lm:omegacb}, 
we have 
\begin{align}
\psi_* (\omega_{X | B} ) 
&\simeq 
\mathcal{H}\mathrm{om}_{\mathcal{O}_{\mathcal{C}}} 
(\psi_* \mathcal{O}_X ,\, \omega_{\mathcal{C} | B}) \notag \\
&\simeq 
\mathcal{H}\mathrm{om}_{\mathcal{O}_{\mathcal{C}}} 
( \mathcal{O}_{\mathcal{C}} ,\, \omega_{\mathcal{C} | B})
\oplus
\mathcal{H}\mathrm{om}_{\mathcal{O}_{\mathcal{C}}} 
(\mathcal{O}_{\mathcal{C}} (-2) \otimes 
\pi_{\mathcal{C}}^* (V_2^{-}),\, 
\omega_{\mathcal{C} | B}) . \notag
\end{align}

So a local section to $\psi_* (\omega_{X | B} )$ is specified 
by a pair $(s, t)$, where $s$ is a local section to 
$\mathcal{O}_{\mathcal{C}}(-1) \otimes \pi_{\mathcal{C}}^* (\mathcal{M})$
identified with that to  
$\mathcal{H}\mathrm{om}_{\mathcal{O}_{\mathcal{C}}} 
( \mathcal{O}_{\mathcal{C}} ,\, \omega_{\mathcal{C} | B})$, 
and $t$ is a local section to 
$\mathcal{O}_{\mathcal{C}}(1) \otimes
\pi_{\mathcal{C}}^* (\mathcal{M} \otimes (V_2^{-})^{\otimes (-1)})$
identified with that to 
$\mathcal{H}\mathrm{om}_{\mathcal{O}_{\mathcal{C}}} 
(\mathcal{O}_{\mathcal{C}} (-2) \otimes 
\pi_{\mathcal{C}}^* (V_2^{-}),\, 
\omega_{\mathcal{C | B}}) $. 
In the same way, a local section to 
$\psi_* \mathcal{O}_X \simeq 
\mathcal{O}_{\mathcal{C}} \oplus 
(\mathcal{O}_{\mathcal{C}} (-2) \otimes \pi_{\mathcal{C}}^* (V_2^{-}))$
is specified by a pair $(\alpha , \beta)$, where 
$\alpha$ is a local section to $\mathcal{O}_{\mathcal{C}}$, and 
$\beta$ is a local section to 
$\mathcal{O}_{\mathcal{C}} (-2) \otimes \pi_{\mathcal{C}}^* (V_2^{-})$.

One can easily check that via the identifications above the 
$\psi_* \mathcal{O}_X$--module structure of 
$\psi_* (\omega_{X | B} )$ is given locally by 
\[
(\alpha, \beta) (s, t) = (\alpha s + \beta t, \, \alpha t + \beta \delta s),
\] 
where 
\[
\delta \in \mathrm{Hom} ((V_2^-)^{\otimes 2},\, \mathcal{A}_4)
\simeq H^0 (
\mathcal{A}_4 \otimes (V_2^-)^{\otimes (-2)}) 
\simeq H^0 (
\mathcal{O}_{\mathcal{C}} (4) 
\otimes \pi_{\mathcal{C}}^* ((V_2^-)^{\otimes (-2)})) 
\]
is the 
global section to 
$\mathcal{H}\mathrm{om}_{\mathcal{O}_B} ((V_2^-)^{\otimes 2},\, \mathcal{A}_4)$
(identified with that to 
$\mathcal{O}_{\mathcal{C}} (4) 
\otimes \pi_{\mathcal{C}}^* ((V_2^-)^{\otimes (-2)})$) 
introduced just after Proposition \ref{prop:amodstrofr}. 

From this it follows that $\psi_* (\omega_{X | B} )$ is locally a 
free $\psi_* \mathcal{O}_X$--module of rank $1$, and we can take 
$(0, u)$ as its local base, with $u$ being a local base 
of the invertible $\mathcal{O}_{\mathcal{C}}$--module 
$\mathcal{O}_{\mathcal{C}}(1) \otimes 
\pi_{\mathcal{C}}^* (\mathcal{M} \otimes (V_2^{-})^{\otimes (-1)})$. 
If we use this to compute the natural morphism 
$\psi^* (\mathcal{O}_{\mathcal{C}}(1) \otimes 
\pi_{\mathcal{C}}^* (\mathcal{M} \otimes (V_2^{-})^{\otimes (-1)})) 
\to \omega_{X | B}$  
induced by  (\ref{eqn:incltodrctim}), we see easily that 
this induced morphism  
$\psi^* (\mathcal{O}_{\mathcal{C}}(1) \otimes 
\pi_{\mathcal{C}}^* (\mathcal{M} \otimes (V_2^{-})^{\otimes (-1)})) 
\to \omega_{X | B}$ is an isomorphism. 
Hence we have the assertion. \qed    
 
\begin{corollary}   \label{cor:omeganformula}
For any integer $n \geq 0$, 
\begin{multline}
\psi_* (\omega_{X | B}^{\otimes n}) \simeq    
(\mathcal{O}_{\mathcal{C}}(n) \otimes 
\pi_{\mathcal{C}}^* (\mathcal{M} \otimes (V_2^{-})^{\otimes (-1)})^{\otimes n})
\oplus \\
(\mathcal{O}_{\mathcal{C}}(n-2) \otimes 
 \pi_{\mathcal{C}}^* (V_2^-)  \otimes 
\pi_{\mathcal{C}}^* (\mathcal{M} \otimes (V_2^{-})^{\otimes (-1)})^{\otimes n})
\notag
\end{multline}
holds. Moreover, in the above, 
\begin{align}
&\psi_* (\omega_{X | B}^{\otimes n})^{\pm}  
\simeq
\mathcal{O}_{\mathcal{C}}(n) \otimes 
\pi_{\mathcal{C}}^* (\mathcal{M} \otimes (V_2^{-})^{\otimes (-1)})^{\otimes n}
\label{eqn:omnplusminus}   \\
&\psi_* (\omega_{X | B}^{\otimes n})^{\mp} 
\simeq 
\mathcal{O}_{\mathcal{C}}(n-2) \otimes 
 \pi_{\mathcal{C}}^* (V_2^-)  \otimes 
\pi_{\mathcal{C}}^* (\mathcal{M} \otimes (V_2^{-})^{\otimes (-1)})^{\otimes n}
\label{eqn:omnminusplus}
\end{align} 
hold, where 
$\psi_* (\omega_{X | B}^{\otimes n}) 
= 
\psi_* (\omega_{X | B}^{\otimes n})^+ 
\oplus
\psi_* (\omega_{X | B}^{\otimes n})^-$ 
is the decomposition into eigen-sheaves with respect to the 
action by the hyperelliptic involution of $\bar{f}$.
\end{corollary}

Proof. 
This follows from 
$\psi_* \mathcal{O}_X \simeq 
\mathcal{O}_{\mathcal{C}} \oplus 
(\mathcal{O}_{\mathcal{C}} (-2) \otimes \pi_{\mathcal{C}}^* (V_2^{-}))$, 
Lemma \ref{lm:omegaxb}, and the projection formula.  \qed 

\begin{proposition}   \label{prop:v2-descr}
Let $V_2 = V_2^+ \oplus V_2^-$ be the decomposition as in the beginning of 
this section, and $L$, the kernel of 
$\sigma_2 : \mathrm{Sym}^2 V_1 \to V_2^+$.
Then 
\[
V_2^- \simeq (\det V_1) \otimes L^{\otimes (-1)} .
\]
\end{proposition}

Proof. 
Note that since $X$ has at most rational double points as 
its singularities, we have 
$V_n = f_* (\omega_{S | B}^{\otimes n}) 
\simeq \bar{f}_* (\omega_{X | B}^{\otimes n})$
for any integer $n \geq 0$. 
Thus, by taking $\bar{f}_*$ of (\ref{eqn:omnplusminus}) and 
(\ref{eqn:omnminusplus}) of Corollary \ref{cor:omeganformula}, 
we obtain
\begin{align}
V_1 &= \bar{f}_* (\omega_{X | B}) 
      \simeq \mathcal{A}_1 
            \otimes (\mathcal{M} \otimes (V_2^-)^{\otimes (-1)}) 
              \label{eqn:v1} \\
V_2^- &\simeq V_2^- \otimes   
 (\mathcal{M} \otimes (V_2^{-})^{\otimes (-1)})^{\otimes 2} 
              \label{eqn:v2-}
\end{align}

Since the locally free sheaf $V_1 \simeq \mathcal{A}_1$ has 
rank $3$, by taking the determinant bundles 
of both hands of (\ref{eqn:v1}), 
we obtain 
$(\mathcal{M} \otimes (V_2^{-})^{\otimes (-1)})^{\otimes 3}
\simeq \mathcal{O}_B$.  And also, since $V_2^-$ is invertible, 
by (\ref{eqn:v2-}), we obtain 
$(\mathcal{M} \otimes (V_2^{-})^{\otimes (-1)})^{\otimes 2}
\simeq \mathcal{O}_B$. 
Thus we obtain 
$\mathcal{M} \otimes (V_2^{-})^{\otimes (-1)}
\simeq \mathcal{O}_B$. 
Since $\mathcal{M} = (\det V_1) \otimes L^{\otimes (-1)}$, 
this implies the assertion.    \qed

\section{Proof of the main theorem}   \label{scn:pfmaintheorem}

In this section, we shall prove our main theorem, i.e., 
Theorem \ref{thm:maintheorem}. 
Let us begin with the proof of 
Propositions \ref{prop:5tuplefibr} and \ref{prop:fibr5tuple}.
\medskip

{\sc Proof of Proposition \ref{prop:fibr5tuple}}
\smallskip

Let $f: S \to B$ be a hyperelliptic fibration of genus $3$ 
as in Proposition \ref{prop:fibr5tuple}, and 
$\bar{f} : X \to B$, its relative canonical model. 
Then $X$ has at most rational double points as its singularities. 
Thus by the same argument 
as in Catanese--Pignatelli \cite[Theorem 4.7]{fibrationsI'},
$\mathcal{C} = \mathcal{P}\mathrm{roj}\, \mathcal{A}$ also 
has at most rational double points as its singularities. 
From this together with Propositions \ref{prop:amodstrofr} 
and \ref{prop:v2-descr}, we see that 
$(B, V_1, V_2^+, \sigma_2, \delta)$ for our $f$ is an 
admissible $5$-tuple, hence the assertion.   \qed 
\medskip

{\sc Proof of Proposition \ref{prop:5tuplefibr}}
\smallskip

Let $(B, V_1, V_2^+, \sigma_2, \delta)$ be an admissible 
$5$-tuple. Construct for this $5$-tuple the commutative diagram 
\begin{equation*} 
\begin{CD}  
 S  @>\text{$\pi$}>>  X  @>\text{$\psi$}>> \mathcal{C}  @>\text{$$}>> 
   \mathbb{P} \\
 @VV\text{$f$}V  @VV\text{$\bar{f}$}V  @VV\text{$\pi_{\mathcal{C}}$}V @VV\text{$\pi_{\mathbb{P}}$}V\\
   B    @>\text{$=$}>> B   @>\text{$=$}>>  B @>\text{$=$}>> B \, , 
\end{CD}
\end{equation*}
following the procedure in Subsection \ref{subscn:5tofibr}. 
By repeating the computations in Section \ref{scn:strrelcanalg}, 
we see easily that this $f: S \to B$ is a hyperelliptic fibration  
of genus $3$ and that 
$\omega_{S | B} \simeq (\psi \circ \pi)^* \mathcal{O}_{\mathcal{C}} (1)$.
This computation of $\omega_{S | B}$ implies that 
the sheaf $\omega_{S | B} |_F  = \mathcal{O}_F (K_F)$ 
is generated by global sections 
for any fiber $F$ of $f: S \to B$. 
So the fibration $f: S \to B$ is relatively minimal.   

Thus we only need to show that all the fibers of $f: S \to B$ 
are $2$-connected. But this follows from the following lemma: 

\begin{lemma}   \label{lm:2connequivcond}
Let $F$ be a fiber of a relatively minimal hyperelliptic 
fibration of genus $3$. Then the natural morphism 
\[
\phi : 
Sym^2 H^0 (\mathcal{O}_F (K_F)) \to H^0 (\mathcal{O}_F (2 K_F))^+ 
\]
is non surjective, if and only if $F$ fails to be $2$-connected. 
\end{lemma}

This lemma follows from the computation of the canonical rings 
of genus $3$ fibers in the Thesis of Mendes Lopes \cite{thesismlopes}. 
As shown by Konno and Mendes Lopes in \cite[Example $1$]{questionofreid} 
(see also \cite[Theorem III]{1-2-3}), 
the computation in \cite{thesismlopes} for multiple fibers 
includes false results. 
This however does not affect our lemma. To avoid to force 
readers to check and patch the pieces from proof in \cite{thesismlopes}, 
we give here a proof of Lemma \ref{lm:2connequivcond}  
in such a way that the correctness of the lemma is clear.  

Proof of Lemma \ref{lm:2connequivcond}. 
We have already seen that the $2$-connectedness of $F$ implies 
the surjectivity of the morphism in the assertion
(see the paragraph just after Lemma \ref{lm:baseofrfkf}). 
Thus we only need to show the $2$-connectedness of $F$ 
assuming the surjectivity of the morphism.  
So assume contrary that the fiber $F$ is non $2$-connected. 
We shall show the non-surjectivity of the morphism for this case.     

If $F$ is $1$-connected (but non $2$-connected), 
then by \cite[Theorem 1.18, Chapter III]{thesismlopes}, 
the image of the multiplication map 
$\mathrm{Sym}^2 H^0 (\mathcal{O}_F (K_F)) 
\to H^0 (\mathcal{O}_F (2 K_F))$ has codimension $\geq 2$ in 
$H^0 (\mathcal{O}_F (2 K_F))$.  Thus $\phi$ is non surjective 
in this case. 

So assume that $F$ is non $1$-connected. 
Then there exists a $1$-connected divisor $D$ with 
$p_a (D) = 2$ such that $F = 2D$. 
By the standard short exact sequence 
\[
0 \to \mathcal{O}_D (nK_F - D) \to \mathcal{O}_F (nK_F) 
            \to \mathcal{O}_D (nK_F) \to 0,       
\]
we see that $H^0 (\mathcal{O}_F (nK_F)) \to H^0 (\mathcal{O}_D (nK_F))$ 
is surjective for any $n \geq 1$ 
(for the proof for the case $n=1$, use the standard fact 
that $\mathcal{O}_D (D)$ is 
a non-trivial $2$-torsion in $\mathrm{Pic}^0 (D)$). 
By the Riemann--Roch theorem, we see moreover 
$h^0 (\mathcal{O}_D (K_F)) = 1$ and $h^0 (\mathcal{O}_D (2K_F)) = 3$.

Now consider the natural commutative diagram 
\begin{equation} \notag
\begin{CD}  
   H^0 (\mathcal{O}_F (K_F))^{\otimes 2}    @>\text{$$}>> H^0 (\mathcal{O}_D (K_F))^{\otimes 2}\\
   @VV\text{$\lambda$}V  @VV\text{$$}V\\
  H^0 (\mathcal{O}_F (2 K_F))  @>\text{$\mu$}>> H^0 (\mathcal{O}_D (2 K_F)), 
\end{CD}
\end{equation}
where the vertical arrows are the multiplication maps, and the 
horizontal arrows are the restriction maps. 
The image $\mathrm{Im}\, ( \mu \circ \lambda )$ has codimension 
$\geq 2$ in $H^0 (\mathcal{O}_D (2 K_F))$, 
since the composite $\mu \circ \lambda$ factors through 
$H^0 (\mathcal{O}_D (K_F))^{\otimes 2} \to H^0 (\mathcal{O}_D (2 K_F))$ 
(note that we have $h^0 (\mathcal{O}_D (K_F)) = 1$).  
Since $\mu$ is surjective, this implies that 
the image $\mathrm{Im}\, \lambda$ has codimension $\geq 2$
in $H^0 (\mathcal{O}_F (2 K_F))$. 
Thus the image of $\phi$ has codimension $\geq 1$ 
in $H^0 (\mathcal{O}_F (2 K_F))^+$, i.e., 
$\phi$ is non surjective in this case, too. 

Summing these up together, we see 
that the fiber $F$ needs to be $2$-connected, 
if $\phi$ is surjective. \qed
\medskip

Now let us prove Theorem \ref{thm:maintheorem}.
\medskip

{\sc Proof of Theorem \ref{thm:maintheorem}}
\smallskip

By the computations in Section \ref{scn:strrelcanalg}, 
we see the existence of the one--to--one correspondence 
as in the assertion; in fact, given a fibration $f: S \to B$ 
as in the assertion, Proposition \ref{prop:amodstrofr} 
and Remark \ref{rem:multstrdscdbydelta} 
implies that the fibration $f: S \to B$ is completely recovered by 
the associated $5$-tuple, and by repeating the argument in 
Section \ref{scn:strrelcanalg} 
we see that the association in Proposition \ref{prop:5tuplefibr} 
and that in Proposition \ref{prop:fibr5tuple} are mutually inverse.  

Thus we only need to show the assertion concerning the numerical 
invariants of the surface $S$. But this follows from the 
following standard formula for $\deg V_n$:
\begin{align}
&\chi (\mathcal{O}_S) - (3-1)(b-1) = \deg V_1 \label{eqn:chiformula}\\
&c_1^2 (S) - 8(3-1)(b-1) = \deg V_2 - \deg V_1, \label{eqn:c1sqformula}
\end{align}
where $b$ is the genus of the base curve $B$. 
In fact, from the equality $\deg V_2 = \deg V_2^+ + \deg V_2^-$,
the short exact sequence 
$0 \to L \to \mathrm{Sym}^2 V_1 \to V_2^+ \to 0$, 
and the isomorphism 
$V_2^- \simeq (\det V_1) \otimes L^{\otimes (-1)}$, 
we infer 
\[
\deg V_2^+ = 4 \deg V_1 - \deg L \qquad \deg V_2^- = \deg V_1 - \deg L,
\]
hence $\deg V_2 = 5 \deg V_1 - 2 \deg L$. 
By this together with (\ref{eqn:chiformula}) and (\ref{eqn:c1sqformula}), 
we obtain the assertion. \qed

\section{Sufficient conditions for the existence of admissible $5$-tuples and some examples} \label{scn:existenceandexs}

In this section, we shall study conditions for the existence of 
admissible $5$-tuples for the case $B \simeq \mathbb{P}^1$,  
show for some cases the uniqueness of the deformation type and 
the simply connectedness of the resulting 
surfaces $S$'s, and 
give some examples of our fibrations $f: S \to B$'s    
including those with minimal regular $S$'s with 
$p_g = 4$ and $c_1^2 = 8$. 
This last examples with $p_g = 4$ and $c_1^2 = 8$ turn out 
to belong to the family $\mathcal{M}_0$ 
in Bauer--Pignatelli \cite{caninvpg4c8}. 

In this section, we assume $B \simeq \mathbb{P}^1$. 
Then by Grothendieck's theorem and Fujita's semi-positivity 
theorem, there exist integers 
$0 \leq d_0 \leq d_1 \leq d_2$ and $e_0$ such that 
$V_1 = \bigoplus_{\lambda = 0}^2 \mathcal{O}_B (d_{\lambda})$, 
$L = \mathcal{O}_B (e_0)$.    
In what follows for a relatively minimal genus $3$ hyperelliptic 
fibration $f: S \to B=\mathbb{P}^1$ as in our main theorem, we let 
$\chi_f  = \chi (\mathcal{O}_S) + 2$ and $K_f^2 = K_S^2 + 16$. 
By our main theorem, we have 
\[
\chi_f = \deg V_1 
\qquad K_f^2 = 4 \deg V_1 - 2 \deg L. 
\] 

First note that if we denote by $l$ and $x_{\lambda}$ local bases 
of the invertible sheaves $L$ and $\mathcal{O}_B (d_{\lambda})$ 
($0 \leq \lambda \leq 2$), respectively, then any 
$\mathcal{O}_B$-module homomorphism  
$L = \mathcal{O}_B (e_0) \to \mathrm{Sym}^2 V_1
= \bigoplus_{i_0 + i_1 + i_2 = 2} 
\mathcal{O}_B (\sum_{\lambda = 0}^2 i_{\lambda} d_{\lambda}) $ 
is given by 
\begin{equation} \label{eqn:exprltosym2v1}
\varPhi_{(\alpha_{i_0\, i_1\, i_2})} : l \mapsto 
(\sum_{i_0 + i_1 + i_2 = 2} \alpha_{i_0 i_1 i_2}) \, l
= \sum_{i_0 + i_1 + i_2 = 2} a_{i_0 i_1 i_2} x_0^{i_0}x_1^{i_1}x_2^{i_2},
\end{equation}
for global sections $\alpha_{i_0 i_1 i_2} 
\in H^0 (\mathcal{O}_B (\sum_{\lambda =0}^2 i_{\lambda} d_\lambda - e_0))$ 
(note that $\alpha_{i_0 i_1 i_2} l 
= a_{i_0 i_1 i_2} x_0^{i_0}x_1^{i_1}x_2^{i_2}$ is a local section 
to $\mathcal{O}_B (\sum_{\lambda =0}^2 i_{\lambda} d_\lambda)$). 
In our case we have $\deg V_1 = \sum_{\lambda = 0}^2 d_{\lambda}$. 
So we restrict our argument in this section to the case 
$\sum_{\lambda = 0}^2 d_{\lambda} > 0$, i.e., the case of 
fibrations with non-constant moduli. 

\begin{proposition} \label{prop:necescond}
Let $0 \leq d_0 \leq d_1 \leq d_2$ and $e_0$ be integers such that 
$\sum_{\lambda = 0}^2 d_{\lambda} > 0$, and assume that 
$B \simeq \mathbb{P}^1$. 
If there exists an admissible $5$-tuple such that 
$V_1 \simeq \bigoplus_{\lambda = 0}^2 \mathcal{O}_B (d_{\lambda})$ 
and $L \simeq \mathcal{O}_B (e_0)$, then 
$e_0 \leq \min \{ d_0 + d_2, \, 2 d_1 \}$
\end{proposition}

Proof.  If we write our $L \to \mathrm{Sym}^2 V_1$ of a 
$5$-tuple as (\ref{eqn:exprltosym2v1}), 
the defining equation of 
$\mathcal{C} \in |\mathcal{O}_{\mathbb{P}}(2) 
\otimes \pi_{\mathbb{P}}^* (L^{\otimes (-1)}) |$ 
is given locally by 
\[
\frac{\sum_{i_0 + i_1 + i_2 = 2} a_{i_0 i_1 i_2} x_0^{i_0}x_1^{i_1}x_2^{i_2}}{1}
\otimes l^{\otimes (-1)} = 0. 
\]
Since we have $\alpha_{i_0 i_1 i_2} 
\in H^0 (\mathcal{O}_B (\sum_{\lambda =0}^2 i_{\lambda} d_\lambda - e_0))$,
if $d_0 + d_2 < e_0$, the global sections 
$\alpha_{2\, 0\, 0}$, $\alpha_{1\, 1\, 0}$, and $\alpha_{1\, 0\, 1}$ 
are identically zero. Then the corresponding $\mathcal{C}$ has 
$1$-dimensional singular locus and this $5$-tuple is not 
admissible.  So we obtain $e_0 \leq d_0 + d_2$ for an admissible $5$-tuple. 
In the same way, we obtain $e_0 \leq 2 d_1$ for an admissible $5$-tuple, too.
\qed

Note that we have
\[
3 K_f^2 - 8 \chi_f = 2 (2(d_0 + d_2 - e_0) + (2d_1 - e_0)).
\]
Thus by Proposition \ref{prop:necescond}, we obtain the 
well known inequality 
$3 K_f^2 - 8 \chi_f \geq 0$ (but only under our restrictive 
assumptions). 

\begin{remark} \label{rm:lowestslope}
By the above, if $3 K_f^2 - 8 \chi_f = 0$, then 
we have $d_0 + d_2 = e_0$ and $2d_1 = e_0$.
Thus there exist integers $d >0$ and $m$ 
such that 
\[
V_1  = \mathcal{O}_B (d-m) 
       \oplus \mathcal{O}_B (d) 
       \oplus\mathcal{O}_B (d+m) \qquad L = \mathcal{O}_B (2d).
\]
Assume $m > 0$. Then the argument above shows that 
the equation of $\mathcal{C}$ in $\mathbb{P}$ is 
$(a_{1\, 0\, 1} x_0 x_2 + a_{0\, 2\, 0} x_1^2 + a_{0\, 1\, 1} x_1 x_2 
+ a_{0\, 0\, 2} x_2^2) \otimes l^{\otimes (-1)} = 0$, 
where $a_{1\, 0\, 1}$ and $a_{0\, 2\, 0}$ are nowhere vanishing 
functions (because $d_0 + d_2 - e_0 = 2d_1 - e_0 = 0$). 
In particular, the section $\varDelta = \{ x_1 = x_2 = 0 \} \subset \mathbb{P}$ 
of $\pi_{\mathbb{P}} : \mathbb{P} \to B$ is contained in our 
relative conic $\mathcal{C}$. 
By an easy computation we see that $\mathcal{C} \to B = \mathbb{P}^1$ is 
the Hirzebruch--Segre surface $\varSigma_m$ with minimal section 
$\varDelta_0 = \varDelta |_{\mathcal{C}}$ and with $\varDelta_0^2 = -m$. 
Let us denote by $\varGamma$ a fiber of $\mathcal{C} \to B$.
Since 
$\omega_{\mathcal{C}} \simeq (\mathcal{O}_{\mathbb{P}} (-1) \otimes 
\pi_{\mathbb{P}}^* ((\det V_1) \otimes L^{\otimes (-1)} \otimes \omega_B)) 
|_{\mathcal{C}}
\simeq \mathcal{O_C} (-2 \varDelta_0 - (2+m) \varGamma )$, 
we obtain $\mathcal{O}_{\mathbb{P}}(1) |_{\mathcal{C}} \simeq 
\mathcal{O}_C (2 \varDelta_0 + (m+d) \varGamma )$. 
Since the branch divisor $D_0$ of the projection 
$S \to \mathcal{C} \simeq \varSigma_m$ belongs to the 
linear system $|\mathcal{O}_{\mathbb{P}}(4) 
\otimes \pi_{\mathbb{P}}^* ((V_2^-)^{\otimes (-2)})| |_{\mathcal{C}}$, 
we see that our surface $S$ is the minimal desingularization 
of the double cover of $\mathcal{C} \simeq \varSigma_m$ 
($0 \leq m \leq 2d/3$) branched along a divisor 
$D_0 \in |8 \varDelta_0 + 2 (2m+d) \varGamma|$ with 
at most negligible singularities. The final description of $S$ is 
valid also for the case $m =0$.  
\end{remark}

\begin{remark} 
In the description of surfaces with the lowest slope in the remark above, 
one can easily check that if $d \leq 2$ then the surface $S$ has 
Kodaira dimension at most $1$ except for the case $d= 2$ and $m=1$. 
If $d=2$ and $m=1$, then $S$ is not minimal, and the minimal model 
$S^*$ has $K_{S^*}^2 = 2$ and $\chi (\mathcal{O}_{S^*}) =4$, hence 
$K_{S^*}^2 = 2 \chi (\mathcal{O}_{S^*}) - 6$. Meanwhile, if $d \geq 3$, 
then $S$ is of general type. For example if $d =3$, then we have 
$m \leq 2$.  In this case, $S$ is a minimal surface with 
$K_S^2 = 8$ and $\chi (\mathcal{O}_S) = 7$
(hence again $K_S^2 = 2 \chi (\mathcal{O}_S) - 6$) 
except for the case $m= 2$. If $d=3$ and $m=2$, then $S$ is not minimal, 
and the minimal model $S^*$ has 
$K_{S^*}^2 = 9$ and $\chi (\mathcal{O}_{S^*}) =7$, hence 
$K_{S^*}^2 = 2 \chi (\mathcal{O}_{S^*}) - 5$. 
It easy to check the following: 
if $d=2$ and $m=1$, $S^*$ is in Case i), 
Theorem 1.6. of \cite{smallc1-1}; 
if $d=3$ and $m=0$, $S$ is in Case iii), 
Theorem 1.6 of \cite{smallc1-1};
if $d=3$ and $m=1$, $S$ is in Case ii), 
Theorem 1.6 of \cite{smallc1-1};
if $d=3$ and $m=2$, $S^*$ is in Case B.2, 
Theorem 1.3 of \cite{smallc1-2}. 
In Case $d=3$ and $m=1$, Horikawa's description says 
our surface $S$ is a minimal resolution of the double cover 
of $\mathbb{P}^2$ branched along a curve of degree $10$ having  
at most negligible singularities.  Our case of genus $3$ fibration 
corresponds to the case where this degree $10$ curve has at least  
one negligible singularity.  Blowing up $\mathbb{P}^2$ at this 
singularity we obtain a projection from the Hirzebruch--Segre surface 
$\varSigma_1 \to \mathbb{P}^2$. The projection $S \to \mathbb{P}^2$ 
lifts to $S \to \varSigma_1 \simeq \mathcal{C}$, 
and our genus $3$ fibration comes from 
the ruling $\varSigma_1 \to B = \mathbb{P}^1$ 
of the Hirzebruch--Segre surface. 
For all other cases, where the fibration comes from is immediately seen
in Horikawa's description.

\end{remark}

Next, let us study sufficient conditions for the existence of 
admissible $5$-tuples. 
The following Lemma follows easily from direct computations. 

\begin{lemma} \label{lm:relconictilde}
Consider $\mathcal{O}_B$-module homomorphism 
$\tilde{\varPhi}_{(\alpha_{i_0\, i_1\, i_2})} : L \to \mathrm{Sym}^2 V_1$ given by 
\[
l \mapsto 
(\alpha_{1\, 0\, 1} + 
\sum_{i_1 + i_2 = 2} \alpha_{0\, i_1\, i_2} )\,l  
=
(a_{1\, 0\, 1} x_0 x_2 + 
\sum_{i_1 + i_2 = 2} a_{0\, i_1\, i_2} x_1^{i_1} x_2^{i_2}),  
\]
and let $\mathcal{C} \subset \mathbb{P} = \mathbb{P} (V_1)$ 
be the relative conic associated to this 
$L \to \mathrm{Sym}^2 V_1$. Let $p \in B$. Then we have the 
following$:$

1) If $\alpha_{1\, 0\, 1} (p) \neq 0$ and $\alpha_{0\, 2\, 0} (p) \neq 0$, 
then the fiber over $p$ of $\mathcal{C} \to B =\mathbb{P}^1$ is 
non-singular.

2) Assume $\alpha_{0\, 2\, 0} (p) = 0$. If $\alpha_{1\, 0\, 1} (p) \neq 0$, 
then the only possible singularity of $\mathcal{C}$ lying over $p$ is 
$(p,\, (x_0: x_1: x_2)) = (p,\, (-a_{0 \, 1\, 1} (p) : a_{1 \, 0\, 1} (p) : 0))$.
This point is a singularity of $\mathcal{C}$ if and only if 
$n: = \mathrm{ord}_p \alpha_{0 \, 2 \, 0} - 1 \geq 1$. In this last case it 
is a singularity of type $A_n$ of $\mathcal{C}$.    

3) Assume $\alpha_{1\, 0\, 1} (p) = 0$. 
If $\alpha_{0\, 0\, 2} (p) \neq 0$ and 
$(4 \alpha_{0\, 2\, 0} \alpha_{0\, 0\, 2} - (\alpha_{0\, 1\, 1})^2) (p) \neq 0$
are satisfied, and neither $\alpha_{1\, 0\, 1}$ nor $\alpha_{0\, 2\, 0}$ is 
identically zero, then the unique singularity of $\mathcal{C}$ lying 
over $p \in B$ is 
$(p, \, (x_0: x_1: x_2)) = (p, \, (1: 0: 0))$. 
This point is a singularity of type $A_n$ of $\mathcal{C}$, 
where $n : = 2 \mathrm{ord}_p \alpha_{1\, 0\, 1} 
+ \mathrm{ord}_p \alpha_{0\, 2\, 0} - 1$ $($$\geq 1$$)$.    

\end{lemma}

Let us use this lemma to prove the following:

\begin{proposition} \label{prop:existence}
Let $B = \mathbb{P}^1$, and 
let $0 \leq d_0 \leq d_1 \leq d_2$ and $e_0$ be 
integers such that $\sum_{\lambda = 0}^2 d_{\lambda} > 0$. 
Assume either \smallskip

A$)$
 $\sum_{\lambda = 0}^2 d_{\lambda} - 2 d_0 
\leq e_0 \leq \min \{ d_0 + d_2, \, 2 d_1 \}$, or 

B$)$ $\sum_{\lambda = 0}^2 d_{\lambda} - \varepsilon / 2
\leq e_0 \leq 2 d_0$, where 
$\varepsilon  : = \min \{ d_0 + 3 d_2, \, 3d_1 + d_2 \}$. \smallskip

Then there exists an admissible $5$-tuple 
(for relatively minimal $2$-connected genus $3$ hyperelliptic fibrations) 
such that 
$V_1 \simeq \bigoplus_{\lambda = 0}^2 \mathcal{O}_B (d_{\lambda})$ 
and $L \simeq \mathcal{O}_B (e_0)$, where 
$0 \to L \to \mathrm{Sym}^2 V_1 \to V_2^+ \to 0$ is 
the short exact sequence in the structure theorem. 
In case A$)$ an admissible $5$-tuple can be taken in such a 
way that the branch divisor of $S \to \mathcal{C}$ is non-singular. 
In case B$)$ an admissible $5$-tuple can be taken in such a 
way that the variety $\mathcal{C}$ is non-singular.

\end{proposition}

Proof. 
Let us first prove the assertion for Case A). 
Let $d_0$, $d_1$, $d_2$, and $e_0$ be integers as in Case A) of the assertion, 
and put $V_1 = \bigoplus_{\lambda = 0}^2 \mathcal{O}_B (d_{\lambda})$ 
and $L = \mathcal{O}_B (e_0)$. 
For each collection $\alpha_{1\, 0\, 1}$, $\{ \alpha_{0\, i_1\, i_2} \}$
($i_1 + i_2 =2$, $i_1$, $i_2 \geq 0$), let us consider the 
morphism $\tilde{\varPhi}_{(\alpha_{i_0\, i_1\, i_2})} : L \to \mathrm{Sym}^2 V_1$
given in  Lemma \ref{lm:relconictilde}.  
Then under the conditions in A), it is an immediate consequence of 
Lemma \ref{lm:relconictilde} that for general 
$\alpha_{1\, 0\, 1}$ and $\{ \alpha_{0\, i_1\, i_2} \}$ 
the $\mathcal{O}_B$-module 
$V_2^+ = 
\mathrm{Cok} 
(\tilde{\varPhi}_{(\alpha_{i_0\, i_1\, i_2})} : L \to \mathrm{Sym}^2 V_1)$
is locally free, and that the relative conic $\mathcal{C} \subset \mathbb{P}$ 
determined by $\tilde{\varPhi}_{(\alpha_{i_0\, i_1\, i_2})}$ has 
at most rational double points as its singularities
(chose $\alpha_{1\, 0\, 1}$, $\alpha_{0\, 2\, 0}$, 
$\alpha_{0\, 0\, 2}$, and $\alpha_{0\, 1\, 1}$ in this order). 
Thus we only need to show the existence of 
$(V_2^-)^{\otimes 2} \to \mathcal{A}_4$ that induces 
an associated canonical model $X$ with at most rational double points 
as its singularities, where 
$\mathcal{A}_4 = \mathrm{Cok}\, 
(L \otimes \mathrm{Sym}^2 V_1 \to \mathrm{Sym}^4 V_1)$ is the 
invertible sheaf as in the definition of a $5$-tuple.    
For this purpose, it is enough to show that 
for a general $(V_2^-)^{\otimes 2} \to \mathrm{Sym}^4 V_1$ the 
induced branch divisor $\mathscr{D} |_{\mathcal{C}}$ 
of $S \to \mathcal{C}$ is non-singular,  
where $\mathscr{D} \in | \mathcal{O}_{\mathbb{P}}(4) 
\otimes \pi_{\mathbb{P}}^* (V_2^-)^{\otimes (-2)}|$ is the 
divisor of $\mathbb{P}$ determined by the morphism 
$(V_2^-)^{\otimes 2} \to 
\mathrm{Sym}^4 V_1 \simeq {\pi_{\mathbb{P}}}_* \mathcal{O}_{\mathbb{P}} (4)$ 
above. 

Note that we have 
$\mathrm{Sym}^4 V_1 \simeq 
\bigoplus_{j_0 + j_1 + j_2 = 4} 
\mathcal{O}_B (\sum_{\lambda = 0}^2 j_{\lambda} d_{\lambda})$  
and 
$(V_2^-)^{\otimes 2} \simeq 
\mathcal{O}_B (2(\sum_{\lambda = 0}^2 d_{\lambda} - e_0))$. 
Thus any morphism $(V_2^-)^{\otimes 2} \to \mathrm{Sym}^4 V_1$ 
is given as 
\[
\varPsi_{(\beta_{j_0\, j_1\, j_2})} : 
y^{\otimes 2} 
\mapsto (\sum_{j_0 + j_1 + j_2 = 4} \beta_{j_0\, j_1\, j_2}) y^{\otimes 2}  
=  \sum_{j_0 + j_1 + j_2 = 4} b_{j_0\, j_1\, j_2} x_0^{j_0} x_1^{j_1} x_2^{j_2}
\]
for $\beta_{j_0\, j_1\, j_2} 
\in 
H^0 (\mathcal{O}_B (\sum_{\lambda = 0}^2 j_{\lambda} d_{\lambda}) 
\otimes (V_2^-)^{\otimes (-2)})$ 
($j_0 + j_1 + j_2 =4$, $j_0$, $j_1$, $j_2 \geq 0$), 
where $y$ is a local base of $V_2^-$. 
Via the composition with the natural morphism 
$\mathrm{Sym}^4 V_1 \to {\pi_{\mathbb{P}}}_* \mathcal{O}_{\mathbb{P}} (4)$, 
the morphism $\varPsi_{(\beta_{j_0\, j_1\, j_2})}$ 
determines an element in 
$H^0 (\mathcal{O}_{\mathbb{P}} (4) 
\otimes \pi_{\mathbb{P}}^* (V_2^-)^{\otimes (-2)})$.    
We define the morphism
\[
\varPsi : \bigoplus_{j_0 + j_1 + j_2 =4} 
H^0 (\mathcal{O}_B (\sum_{\lambda = 0}^2 j_{\lambda} d_{\lambda}) 
\otimes (V_2^-)^{\otimes (-2)})
\to
H^0 (\mathcal{O}_{\mathbb{P}} (4) 
\otimes \pi_{\mathbb{P}}^* (V_2^-)^{\otimes (-2)})
\]
by 
$(\beta_{j_0\, j_1\, j_2}) \mapsto \varPsi_{(\beta_{j_0\, j_1\, j_2})}$. 

By the condition $\sum_{\lambda = 0}^2 d_{\lambda} - 2 d_0 \leq e_0$, 
we see that the linear system 
$|\mathcal{O}_B (\sum_{\lambda = 0}^2 j_{\lambda} d_{\lambda}) 
\otimes (V_2^-)^{\otimes (-2)}|$ is free from base points 
for any $j_0$, $j_1$, $j_2 \geq 0$ such that 
$\sum_{\lambda = 0}^2 j_{\lambda} = 4$. 
Thus the image $\mathrm{Im}\, \varPsi$ determines a linear system 
free from base points. It follows that, for a general 
member $\mathscr{D}$ of  $| \mathrm{Im}\, \varPsi |$, 
the restriction $D_0 = \mathscr{D} |_{\mathcal{C}}$ to $\mathcal{C}$ 
is non-singular.   
Thus $\varPsi_{(\beta_{j_0\, j_1\, j_2})} : (V_2^-)^{\otimes 2} 
\to  \mathrm{Sym}^ 4 V_1$ ($\to \mathcal{A}_4$) induces an admissible 
$5$-tuple, hence the assertion for Case A).  

The assertion for Case B) can be proved by the same method, but 
with a bit more computation. 
Here we only observe the difference from the proof for Case A), 
and give an outline of the proof to the extent that the readers 
can verify the detail of the proof by themselves. 

First, note that the variety $\mathbb{P} = \mathbb{P} (V_1)$ 
contains two subvarieties 
$\varDelta = \{ x_1 = x_2 = 0\} \subset \mathbb{P}$ and
$\mathscr{D}_1^{\prime} = \{ x_2 = 0 \} \simeq
\mathbb{P} 
(\mathcal{O}_B (d_0) \oplus \mathcal{O}_B (d_1)) \subset \mathbb{P}$. 
The variety $\varDelta$ is a section of the projection 
$\pi_{\mathbb{P}} : \mathbb{P} \to B= \mathbb{P}^1$, 
and the variety $\mathscr{D}_1^{\prime} = \{ x_2 = 0 \}$ 
is a divisor of $\mathbb{P}$
and a member of the linear system 
$|\mathcal{O}_{\mathbb{P}} (1) \otimes \pi_{\mathbb{P}}^* \mathcal{O}_B ( -d_2)|$. 
In case B), instead of $\tilde{\varPhi}_{(\alpha_{i_0\, i_1\, i_2})}$, 
we employ the morphism 
$\varPhi_{(\alpha_{i_0\, i_1\, i_2})} : L \to \mathrm{Sym}^2 V_1$ 
given in (\ref{eqn:exprltosym2v1}).
We define the linear map 
$\varPhi : \bigoplus_{i_0 + i_1 + i_2 = 2} 
H^0 (\mathcal{O}_B (\sum_{\lambda = 0}^2 i_{\lambda} d_{\lambda} - e_0))
\to 
H^0 (\mathcal{O}_{\mathbb{P}}(2) \otimes \pi_{\mathbb{P}}^* L^{\otimes (-1)})$
by $(\alpha_{i_0\, i_1\, i_2}) \mapsto \varPhi_{(\alpha_{i_0\, i_1\, i_2})}$. 
Then by the condition $e_0 \leq 2 d_0$ we see that 
the image $\mathrm{Im}\, \varPhi$ determines a base point free 
linear system.  Thus for a general $(\alpha_{i_0\, i_1\, i_2})$ 
the induced relative conic $\mathcal{C} \subset \mathbb{P}$ 
is non-singular, does not contain $\varDelta$, and intersects 
transversally at (if any) each point of $\mathcal{C} \cap \varDelta$.  

Meanwhile, the control of $(\beta_{j_0\, j_1\, j_2})$ becomes 
a bit more complicated compared to that in Case A). 
In Case B), at the worst case of 
$(d_0, d_1, d_2)$ and $e_0$, the seven sections 
$\beta_{j_0\, j_1\, 0}$ ($j_0 + j_1 =4$), 
$\beta_{3\, 0\, 1}$, and $\beta_{2\, 0\, 2}$ vanish. 
So, in Case B), in stead of $\varPsi_{(\beta_{j_0\, j_1\, j_2})}$, 
we employ the morphism 
$\tilde{\varPsi}_{(\beta_{j_0\, j_1\, j_2})} :
(V_2^-)^{\otimes 2} \to \mathrm{Sym}^4 V_1$ 
\[
y^{\otimes 2} \mapsto 
(\beta_{1\, 0\, 3} + \sum_{j_1 + j_2 = 4, \, j_2 \geq 1} \beta_{0\, j_1\, j_2} )
y^{\otimes 2} \notag 
= 
b_{1\, 0\, 3} x_0 x_2^3 
+ \sum_{j_1 + j_2 = 4, \, j_2 \geq 1} b_{0\, j_1\, j_2} x_1^{j_1} x_2^{j_2}  
\]
associated to $\beta_{1\, 0\, 3} \in 
H^0 (\mathcal{O}_B (d_0 + 3 d_2) 
\otimes (V_2^-)^{\otimes (-2)})$
and 
$\beta_{0\, j_1\, j_2} \in 
H^0 (\mathcal{O}_B (j_1 d_1 + j_2 d_2) 
\otimes (V_2^-)^{\otimes (-2)})$
($j_1 + j_2 = 4$, $j_1 \geq 0$, $j_2 \geq 1$), 
where $y$ is a local base of $V_2^-$ as usual.

Since the equation of the divisor 
$\mathscr{D} \in 
|\mathcal{O}_{\mathbb{P}} (4) \otimes \pi_{\mathbb{P}}^* (V_2^-)^{\otimes (-2)}|$ 
associated to $\tilde{\varPsi}_{(\beta_{j_0\, j_1\, j_2})}$ has $x_2$ 
as a prime factor, we obtain the splitting 
$\mathscr{D} = \mathscr{D}_1^{\prime} + \mathscr{D}_2^{\prime}$, 
where 
$\mathscr{D}_1^{\prime} = \{ x_2 = 0 \} 
\in |\mathcal{O}_{\mathbb{P}} (1) 
\otimes \pi_{\mathbb{P}}^* \mathcal{O}_B (-d_2)|$ 
and 
$\mathscr{D}_2^{\prime}  
\in |\mathcal{O}_{\mathbb{P}} (3) 
\otimes 
\pi_{\mathbb{P}}^* (\mathcal{O}_B (d_2) \otimes (V_2^-)^{\otimes (-2)} )|$.  
So let us denote by 
\[
\varPsi^{\prime}_{(\beta_{j_0\, j_1\, j_2})}: (V_2^-)^{\otimes 2} 
\to (\mathrm{Sym}^3 V_1) \otimes \mathcal{O}_B (d_2)
\simeq {\pi_{\mathbb{P}}}_* 
(\mathcal{O}_{\mathbb{P}} (3) \otimes \pi_{\mathbb{P}}^* \mathcal{O}_B (d_2))
\]
the morphism corresponding to the defining equation in $\mathbb{P}$ 
of the divisor $\mathscr{D}_2^{\prime}$, 
and consider the linear map 
$\varPsi^{\prime} : (\beta_{j_0\, j_1\, j_2}) \mapsto 
\varPsi^{\prime}_{(\beta_{j_0\, j_1\, j_2})} \in 
H^0 (\mathcal{O}_{\mathbb{P}} (3) 
\otimes \pi_{\mathbb{P}}^* (\mathcal{O}_B (d_2) \otimes (V_2^-)^{\otimes (-2)} ))$. 

From the condition $\sum_{\lambda = 0}^2 d_{\lambda} - \varepsilon / 2 \leq e_0$, 
we see the following:   
\smallskip

a) The base locus of the linear system $|\mathrm{Im}\, \varPsi^{\prime}|$ 
determined by the image $\mathrm{Im}\, \varPsi^{\prime}$ is contained in 
the subvariety $\varDelta = \{ x_1 = x_2 = 0 \} \subset \mathbb{P}$.
So for a general $\mathscr{D}_2^{\prime}$, its restriction
$\mathscr{D}_2^{\prime} |_{\mathcal{C}}$ to $\mathcal{C}$ 
is non-singular outside 
$\varDelta \cap \mathcal{C}$.  

b) If we take general $\mathcal{C}$, then the divisor 
$\mathscr{D}_1^{\prime} |_{\mathcal{C}}$ of $\mathcal{C}$ 
is non-singular. 

c) If we take $\mathcal{C}$ sufficiently general, then for a general 
$\mathscr{D}_2^{\prime}$, 
the divisor $(\mathscr{D}_1^{\prime} + \mathscr{D}_2^{\prime}) |_{\mathcal{C}}$
of $\mathcal{C}$ has at each point (if any) of $\varDelta \cap \mathcal{C}$ 
a negligible singularity. 
More precisely, these are singularities of 
type $(x_2/x_0)((x_2/x_0)^2 + (x_1/x_0)^3) = 0$, i.e., 
those corresponding to singularities of type $E_7$ of the double cover
(note here that the condition 
$\sum_{\lambda = 0}^2 d_{\lambda} - \varepsilon / 2 \leq e_0$ ensures 
non-vanishing of general $\beta_{1\, 0\, 3}$ and $\beta_{0\, 3\, 1}$, 
which are coefficients of $(x_2/x_0)^3$ and $(x_1/x_0)^3 (x_2/x_0)$,
respectively, and that $(x_1/x_0)$ and $(x_2/x_0)$ form a system of 
local coordinates of $\mathcal{C}$ around these points).  
 
d) At points outside $\varDelta \cap \mathcal{C}$, general 
$\mathscr{D}_1^{\prime} |_{\mathcal{C}}$ and 
$\mathscr{D}_2^{\prime} |_{\mathcal{C}}$ 
at most intersect each other transversally
(by the defining equation of $\mathscr{D}_2$, we see that,  
for a general $\mathcal{C}$, intersection points of   
$\mathscr{D}_1^{\prime} |_{\mathcal{C}}$ and 
$\mathscr{D}_2^{\prime} |_{\mathcal{C}}$ outside $\varDelta \cap \mathcal{C}$
appear only in the fibers over points $p \in B$'s such that 
$\beta_{0\, 3\, 1} (p) = 0$). 
\smallskip

Thus, for a general $(\beta_{j_0\, j_1\, j_2})$, the divisor 
$\mathscr{D} |_{\mathcal{C}} = 
(\mathscr{D}_1^{\prime} + \mathscr{D}_2^{\prime}) |_{\mathcal{C}}$ 
of $\mathcal{C}$ has at most negligible singularities.  
This implies that $(V_2^-)^{\otimes 2} \to \mathcal{A}_4$ 
associated to a general $(\beta_{j_0\, j_1\, j_2})$ induces 
an admissible $5$-tuple, hence the assertion for Case B). \qed 

\begin{remark}   \label{rem:someexamples}
Note that the cases of the lowest slope $8/3$ are covered by 
Case A) of Proposition \ref{prop:existence}
(see Remark \ref{rm:lowestslope}).  
As is well-known, 
contrary to the case of 
genus $2$ fibrations, the slope of a genus $3$ fibration 
can take higher values even if it has only $2$-connected fibers. 
Using our Proposition \ref{prop:existence}, we can for example show 
that the slope can take the value $s$ for any rational number 
$8/3 \leq s \leq 10/3$ even under the $2$-connected assumption. 
In fact, for a sufficiently large positive integer $d$ such that 
$sd/2$ is an integer, put $d_0 = d_1 = d_2 = d$ and 
$e_0 = (4-s) 3d /2$. 
Then by A) of Proposition \ref{prop:existence} we see the existence 
of an admissible $5$-tuple such that 
$V_1 = \bigoplus_{\lambda =0}^2 \mathcal{O}_B (d_{\lambda})$ and 
$L = \mathcal{O}_B (e_0)$, which yields a $2$-connected genus $3$ 
fibration with slope $s$. 
To construct examples with even higher slope, the condition in B) 
is more useful. Put for example 
$d_0 = d$, $d_1 = 5d$, $d_2 = 7d$, and $e_0 = 2d$
for any positive integer $d$.  Then the existence of an admissible 
$5$-tuple is assured by case B) of Proposition \ref{prop:existence}, 
and the resulting fibration has slope $48/13 = 3.6923 \ldots$.  
\end{remark}

Next, let us show for some cases 
the uniqueness of the deformation 
type and the simply connectedness of 
the resulting surfaces $S$'s.  
As for the following Theorem \ref{thm:simplyconnected}, 
we can show with messy computations 
that the maximal slope covered in this theorem 
is $10/3$, hence $< 3.5$.  So the simply connectedness 
follows from Xiao's results \cite[Lemma 2, Theorem 2]{pi1hyperell}, 
since our fibrations $f: S \to B$'s have no multiple fiber.  
In Proposition \ref{prop:exsistsimplconn} 
(see also Remark \ref{rem:simpconnslope}), 
however, we shall give examples with slope $3.6$. 
In this last case, the simply connectedness does not directly 
follow from \cite[Lemma 2, Theorem 2]{pi1hyperell}. 
To unify the the proof, 
we shall give even in Theorem \ref{thm:simplyconnected}  
a proof of simply connectedness using a result from Catanese 
\cite{catonmodulijdg}, 
which works also for the cases in Proposition 
\ref{prop:exsistsimplconn}.    

\begin{theorem}   \label{thm:simplyconnected} 
Let $f: S \to B = \mathbb{P}^1$ be a relatively minimal hyperelliptic 
fibration of genus $3$ with all fibers $2$-connected. 
Assume 
$f_*(\omega_{S|B}) = V_1 \simeq 
\bigoplus_{\lambda = 0}^2 \mathcal{O}_B (d_{\lambda})$ 
with $\sum_{\lambda = 0}^2 d_{\lambda} > 0$ and 
$d_0 \leq  d_1 \leq d_2$.  
Let $L = \ker \sigma_2$ be the kernel of the morphism 
$\sigma_2 : \mathrm{Sym}^2 \, V_1 \to V_2^+$ as in 
Theorem \ref{thm:maintheorem}. Then if
$\sum_{\lambda = 0}^2 d_{\lambda} - 2 d_0 
\leq e_0 \leq \min \{ d_0 + d_2, 2d_1 \}$, then 
the surface $S$ is topologically simply connected, and  
any two such $S$'s having the same 
$(d_0, \, d_1, \, d_2)$ and $e_0$ are equivalent 
under the deformation of complex structures. 
\end{theorem}

Proof. 
Let us first prove the uniqueness of the 
deformation type. 
Assume that $(d_0, \, d_1, \, d_2)$ and $e_0$ as in 
the assertion are given. 
Fix one $0 \to L \to \mathrm{Sym}^2 V_1 \to V_2^+ \to 0$ 
(exact) such that the associated relative conic 
$\mathcal{C} \subset \mathbb{P}$ has at most rational double points. 
Under this fixed $\mathcal{C}$, deforming the morphism 
$(V_2^-)^{\otimes 2} \to \mathcal{A}_4$ correspond to deforming the 
branch divisor of $S \to \mathcal{C}$.   
Thus with the aid of Tjurina's theorem on simultaneous resolution 
of the family of rational double points we see that under the 
fixed $\mathcal{C}$ the deformation type of the resulting surface 
$S$ does not depend on the choice of admissible 
$(V_2^-)^{\otimes 2} \to \mathcal{A}_4$
(note here that to a general member $D$ of the linear system 
to which the branch divisors of $S \to \mathcal{C}$ belong, 
an admissible $(V_2^-)^{\otimes 2} \to \mathcal{A}_4$ correspond,  
since in our case $D$ has at most negligible singularities).  
We denote by $d(\mathcal{C})$ this deformation type, which depends 
only on the choice of the relative conic $\mathcal{C}$. 

Since $\sum_{\lambda = 0}^2 d_{\lambda} - 2 d_0 \leq e_0$, the 
linear system determined by the image $\mathrm{Im}\, \varPsi$ 
of $\varPsi$ (in the proof of Proposition \ref{prop:existence}) 
is base point free. Thus there exists a member 
$\mathscr{D} \in |\mathrm{Im} \varPsi|$ 
whose restriction $\mathscr{D} |_{\mathcal{C}}$ is non-singular 
and passes no singular point of $\mathcal{C}$.  
Under this fixed $\mathscr{D}$, if $\mathcal{C}^{\prime}$ 
is any sufficiently small deformation of $\mathcal{C}$, 
the deformation type of the surface associated to the pair 
$(\mathcal{C}^{\prime}, \mathscr{D})$ is the same as that of 
the surface associated to the original pair 
$(\mathcal{C}, \mathscr{D})$. 
This means $d(\mathcal{C}) = d(\mathcal{C}^{\prime})$. 
Since the space parametrizing $\mathcal{C}$'s are connected, 
this implies the uniqueness of the deformation type of our $S$.  
 
Now let us prove the simply connectedness of the surface $S$. 
Since we have already proved the uniqueness of the deformation type 
of our $S$, we only need to find an $S$ that is simply connected. 
To do this, let us take $\mathscr{D}$ as above such that 
the restriction $D_0 = \mathscr{D} |_{\mathcal{C}}$ is non-singular, and 
show that the associated surface $S$ is simply connected. 

Let $\tilde{\mathcal{C}} \to \mathcal{C}$ be the minimal 
desingularization of $\mathcal{C}$, and 
$\tilde{D}_0$, the total transform of $D_0$ to $\tilde{\mathcal{C}}$. 
Note that the linear system $|\tilde{D}_0|$ is free from base point
(because $\sum_{\lambda = 0}^2 d_{\lambda} - 2 d_0 \leq e_0$). 
From this together with $\tilde{D}_0^2 = 16 e_0 > 0$,  
we see that the divisor $\tilde{D}_0$ is flexible in the 
sense of Definition 1.4., Catanese \cite{catonmodulijdg}. 
Moreover, our variety $\tilde{\mathcal{C}}$ is rational, hence 
simply connected. Then the simply connectedness of $S$ follows from 
Proposition 1.8., Catanese \cite{catonmodulijdg}. \qed

By a similar argument we can prove the following:

\begin{proposition}  \label{prop:exsistsimplconn}
Let $0 \leq d_0 \leq d_1 \leq d_2$ and $e_0$ be integers 
such that $\sum_{\lambda = 0}^2 d_{\lambda} > 0$. 
Assume either \smallskip

A$)$ $\sum_{\lambda = 0}^2 d_{\lambda} -(3d_0 + d_1)/2 \leq e_0 \leq 2 d_0$, or

B$)$ $\sum_{\lambda = 0}^2 d_{\lambda} - \varepsilon^{\prime} / 2
\leq e_0 = 2 d_0$, 
where $\varepsilon^{\prime} = \min \{ d_0 + 3 d_2, \, 4 d_1 \}$. 
\smallskip

Then there exists a relatively minimal genus $3$ hyperelliptic fibration 
$f : S \to B = \mathbb{P}^1$ with all fibers $2$-connected 
such that 
$V_1 \simeq \bigoplus_{\lambda = 0}^2 \mathcal{O}_B (d_{\lambda})$ 
and $L \simeq \mathcal{O}_B (e_0)$ with 
$S$ topologically simply connected. 
\end{proposition}

Prof. 
Since $3d_0 + d_1 \leq \varepsilon^{\prime} \leq \varepsilon$, 
where $\varepsilon$ is as in Proposition \ref{prop:existence}, 
the existence of a relatively minimal fibration $f: S \to B$ 
with all fibers $2$-connected is assured by 
Proposition \ref{prop:existence}. 
Thus we only need to prove that we can take $f: S \to B$ such
that $S$ is simply connected.  

Before we start, let us note that 
in both Cases A) and B) we have $e_0 > 0$. 
In fact, in Case A), we have 
$e_0 \geq \sum_{\lambda = 0}^2 d_{\lambda} - (3 d_0 + d_1)/2 
= (d_1 - d_0)/ 2 + d_2 \geq 0$, and in Case B) we have 
$2d_0 = e_0 \geq \sum_{\lambda = 0}^2 d_{\lambda} - \varepsilon^{\prime} / 2
\geq \sum_{\lambda = 0}^2 d_{\lambda} - 
((d_0 + 3 d_2) + 4 d_1)/4 = (3 d_0 + d_2 )/4 \geq 0$. 
Thus, in both cases, $e_0 \leq 0$ would imply 
$d_0 = d_1 = d_2 = 0$, which contradicts our assumption 
$\sum_{\lambda = 0}^2 d_{\lambda} > 0$. 

Let us prove the assertion for Case A). The condition $e_0 \leq 2 d_0$ 
ensures that a general $\mathcal{C} \subset \mathbb{P}$ is 
non-singular, and $\mathcal{C}$ intersects 
$\varDelta = \{ x_1 = x_2 =0 \}$ transversally   
at (if any) each point of $\mathcal{C} \cap \varDelta$
(because $\alpha_{2\, 0\, 0} \in H^0 (\mathcal{O}_B (2d_0 - e_0))$, 
$2 d_0 - e_0 \geq 0$). 
Then recall that our $\mathscr{D}$ is locally defined by
\[
(\sum_{j_0 + j_1 + j_2 = 4} \beta_{j_0\, j_1\, j_2}) y^{\otimes 2}  
=  \sum_{j_0 + j_1 + j_2 = 4} b_{j_0\, j_1\, j_2} x_0^{j_0} x_1^{j_1} x_2^{j_2}
=0
\]
as in the proof of Proposition \ref{prop:existence}. 
So the condition $\sum_{\lambda = 0}^2 d_{\lambda} -(3d_0 + d_1)/2 \leq e_0$ 
implies that the image of $\varPsi$ in the proof for case A) of 
Proposition \ref{prop:existence} gives 
a linear system with base locus contained in $\varDelta$
(because of the non-vanishing of general $\beta_{0\, 4\, 0}$ 
and $\beta_{0\, 0\, 4}$), that for a general member $\mathscr{D}$ 
of this linear system its restriction $D_0 = \mathscr{D} |_{\mathcal{C}}$ 
is smooth at these base points (because of the non-vanishing of general  
$\beta_{3\, 0\, 1}$), and that for any two general $\mathscr{D}$'s 
their restrictions $D_0 = \mathscr{D} |_{\mathcal{C}}$'s intersect 
each other transversally at these base points 
(because of the non-vanishing of general  
$\beta_{3\, 1\, 0}$ and $\beta_{3\, 0\, 1}$).   
These together with 
$D_0^2 = 16 e_0 > 0$ ensure that a general 
$D_0 = \mathscr{D} |_{\mathcal{C}}$ is flexible, 
hence the assertion for Case A). 

Let us prove the assertion for case B). 
By the condition $e_0 = 2 d_0$, a general 
$\mathcal{C}$ is non-singular and satisfies 
$\mathcal{C} \cap \varDelta = \emptyset$ 
(because $\alpha_{2\, 0\, 0} \in H^0 (\mathcal{O}_B (2d_0 - e_0))$, 
$2d_0 - e_0 = 0$).  
But, from the condition 
$\sum_{\lambda = 0}^2 d_{\lambda} - \varepsilon^{\prime} / 2
\leq e_0 $, we see again that 
the image $\mathrm{Im}\, \varPsi$ determines a linear system 
with base locus contained in $\varDelta$. 
So the linear system $\{ D_0 = \mathscr{D} |_{\mathcal{C}} \}$ is 
free from base points. 
This together with $D_0^2 = 16 e_0 > 0$ ensures that a general 
$D_0$ is flexible, hence the assertion for Case B). \qed 

\begin{remark}  \label{rem:simpconnslope}
As an example, consider the case  
$d_0 = d$, $d_1 = 4d$, $d_2 = 5d$, and $e_0  = 2d$ for a  
positive integer $d > 0$. This case is covered by 
Case B) of Proposition \ref{prop:exsistsimplconn}. 
Thus we obtain an example of $f: S \to B \simeq \mathbb{P}^1$ 
with $V_1 \simeq \mathcal{O}_B (d) \oplus 
\mathcal{O}_B (4d) \oplus \mathcal{O}_B (5d)$,  
$L \simeq \mathcal{O}_B (2d)$, and topologically simply connected $S$. 
This $S$ has numerical invariants $c_1^2 = 36d -16$ and 
$\chi ({\mathcal{O}_S}) = 10d -2$, hence the slope of $f$ being 
$36/10 = 3.6$, and the ratio $c_1^2 / \chi (\mathcal{O}_S)$ 
asymptotically being $3.6$. 
Note even for the case $c_1^2  < 3 \chi(\mathcal{O}_S)$, 
even the algebraic fundamental group is not necessarily trivial 
(see foe example Mendes Lopes--Pardini \cite{mlopespardini3chi}).    
\end{remark}

\begin{remark}
In Proposition \ref{prop:exsistsimplconn}, 
we put the conditions in order to assure that the 
general divisor $D_0$ is smooth at points in $\varDelta \cap \mathcal{C}$. 
It is obvious that allowing some mild negligible singularities 
at points in $\varDelta \cap \mathcal{C}$, we can weaken the conditions 
in Proposition \ref{prop:exsistsimplconn}. 
We do not execute it here, and instead chose to execute it when we need, 
since the conditions in the results will be a bit more complicated.  
\end{remark}

One of the interesting cases that are covered by our Proposition 
\ref{prop:exsistsimplconn} is the case 
$d_0 = 1$, $d_1 = d_2 =3$, and $e_0 = 2$. 
In this case, for a general admissible $5$-tuple, 
the associated hyperelliptic fibration 
$f: S \to B = \mathbb{P}^1$ has minimal regular S with 
$c_1^2 = 8$ and $p_g = 4$. 
In our case, obviously, the canonical map of our 
$S$ factors through the hyperelliptic involution $\iota$ 
of the fibration $f: S \to B$.  
Thus these $S$ are of one of the types given in the list by 
I.\, Bauer and R.\, Pignatelli \cite{caninvpg4c8} 
of minimal regular surfaces 
with $c_1^2 = 8$ and $p_g = 4$ with canonical involution. 

For the use in a sequel to this paper, 
let us specify of which type in the list in \cite{caninvpg4c8} 
our surfaces $S$'s are.

\begin{proposition} \label{prop:bauerpignatellim0}
Let $d_0 =1$, $d_1 = d_2 =3$, and $e_0 = 2$, 
hence a special case of Case B$)$ of Proposition 
\ref{prop:exsistsimplconn}. 
Let $f : S \to B = \mathbb{P}^1$ be a relatively minimal genus $3$ 
hyperelliptic fibration with all fibers $2$-connected obtained 
by Proposition \ref{prop:exsistsimplconn} using 
a non-singular $\mathcal{C}$ such that 
$\mathcal{C} \cap \varDelta = \emptyset$. 
Then $S$ is a minimal regular surface with 
$c_1^2 = 8$ and $p_g = 4$ whose corresponding point $[S]$ 
in the moduli space lies in the strata $\mathcal{M}_0$ 
in the Main Theorem of Bauer--Pignatelli \cite{caninvpg4c8}. 
\end{proposition}
 
Proof. 
The minimality of $S$ follows from $\mathcal{C} \cap \varDelta = \emptyset$ 
and $f_* (\omega_S) \simeq V_1 \otimes \omega_B \simeq 
\mathcal{O}_B (-1) \oplus \mathcal{O}_B (1) \oplus \mathcal{O}_B (1)$, 
since these imply that the canonical system of $S$ is free from base points. 
The equalities $c_1^2 = 8$, $p_g = 4$, and $q = 0$ follow from 
Theorem \ref{thm:maintheorem} and $p_g = h^0 (V_1 \otimes \omega_B)$. 
Thus we only need to show that the corresponding point 
$[S]$ in the moduli space 
lies in the strata $\mathcal{M}_0$. 

To prove the assertion concerning the corresponding point $[S]$, 
let us first note that, 
in Proposition \ref{prop:exsistsimplconn}, 
if $\tilde{\mathcal{C}} \to \mathcal{C}$ is the minimal 
resolution of singularities of $\mathcal{C}$, there exist  
a Hirzebruch--Segre surface $\bar{\mathcal{C}} \to B$ and 
a projection $\tilde{\mathcal{C}} \to \bar{\mathcal{C}}$ 
compatible with the projections to $B$ 
such that this $\tilde{\mathcal{C}} \to \bar{\mathcal{C}}$ is 
a blowing-up of (possibly infinite near) 
$(2 \sum_{\lambda = 0}^2 d_{\lambda} - 3 e_0)$ points.   
In fact, this follows from 
$K_{\tilde{\mathcal{C}}}^2 
= 8 - (2 \sum_{\lambda = 0}^2 d_{\lambda} - 3 e_0)$, 
which in tern follows from 
$(\mathcal{O}_{\mathbb{P}} (1) |_{\mathcal{C}})^2 
= 2 \sum_{\lambda = 0}^2 d_{\lambda} - e_0$. 
In our case of $d_0 = 1$, $d_1 = d_2 = 3$, and $e_0 = 2$ 
with non-singular $\mathcal{C}$, we have 
$\tilde{\mathcal{C}} = \mathcal{C}$ and the projection 
$\mathcal{C} \to \bar{\mathcal{C}}$ is a blowing-up of 
$8$ points. 
In particular, we obtain $K_C^2 = 0$ for our case.   

Let $\iota$ be the hyperelliptic involution of $f: S \to B$. 
Then the canonical map of $S$ factors through $\iota$, and 
the involution $\iota$ has no isolated fixed point 
because the branch divisor $D_0 = \mathscr{D} |_{\mathcal{C}}$ 
has at most negligible singularities.  
Thus our $S$ corresponds to a point in the strata $\mathcal{M}_0^{\mathrm{div}}$ 
or to a point in the strata $\mathcal{M}_0$ 
(see Section 3 of \cite{caninvpg4c8}, especially, 
Proposition 3.2, Theorem 3.3, and Theorem 3.5). 

Let us show that the corresponding point $[S]$ 
belongs to the strata $\mathcal{M}_0$. 
To do this, since we already have $K_{\mathcal{C}}^2 = 0$, 
we only need to show 
that $\mathcal{C}$ has no $(-1)$-curve contained in the image 
by $S \to \mathcal{C}$ of a fundamental cycle of $S$
(see Proposition 3.2, Theorem 3.5, \cite{caninvpg4c8}).   

Recall we have $V_1 \simeq 
\mathcal{O}_B (1) \oplus \mathcal{O}_B (3) \oplus \mathcal{O}_B (3)$, 
hence 
$f_* (\omega_S) \simeq 
\mathcal{O}_B (-1) \oplus \mathcal{O}_B (1) \oplus \mathcal{O}_B (1)$. 
From this we infer 
$H^0 (\omega_S) \simeq 
H^0 (\mathcal{O}_{\mathbb{P}} (1) 
\otimes \pi_{\mathbb{P}}^* \omega_B) \simeq \mathbb{C}^4$. 
Thus the canonical map of $S$ is the composition of the following 
three maps: the natural projection $S \to \mathcal{C}$, 
the natural inclusion $\mathcal{C} \to \mathbb{P} = \mathbb{P} (V_1)$, 
and the rational map $\varPhi_{|\mathcal{O}_{\mathbb{P}} (1) 
\otimes \pi_{\mathbb{P}}^* \omega_B|} : \mathbb{P} - - \to \mathbb{P}^3$ 
associated to the linear system 
$|\mathcal{O}_{\mathbb{P}} (1) \otimes \pi_{\mathbb{P}}^* \omega_B|$.  
But it easy to see the following: \smallskip

a) the indeterminacy locus of 
$\varPhi_{|\mathcal{O}_{\mathbb{P}} (1) 
\otimes \pi_{\mathbb{P}}^* \omega_B|}$ is 
$\varDelta = \{ x_1 = x_2 = 0 \}$; 

b) the image of $\varPhi_{|\mathcal{O}_{\mathbb{P}} (1) 
\otimes \pi_{\mathbb{P}}^* \omega_B|}$ is a non-singular 
quadric $\mathcal{Q} \simeq B \times \mathbb{P}^1 \subset \mathbb{P}^3$ 
(with $\mathbb{P} - - \to \mathcal{Q}$ 
compatible with the natural projections onto $B$
of $\mathbb{P}$ and $\mathcal{Q} \simeq B \times \mathbb{P}^1$); 

c) for any $p \in B$, the restriction 
$\pi_{\mathbb{P}}^{-1} (p) \simeq \mathbb{P}^2 - - \to 
\{ p \} \times \mathbb{P}^1$ ($\subset B \times \mathbb{P}^1$) 
(to the fibers over $p$) of the rational map 
$\mathbb{P} - - \to \mathcal{Q}$ 
is the linear projection  of 
$\pi_{\mathbb{P}}^{-1} (p) \simeq \mathbb{P}^2$ 
from the center $\pi_{\mathbb{P}}^{-1} (p) \cap \varDelta$. \smallskip

From these together with $\mathcal{C} \cap \varDelta = \emptyset$, 
it follows that $\mathcal{C} \to \mathcal{Q} \subset \mathbb{P}^3$ 
contracts no curve, and is of mapping degree $2$.  
Thus the canonical map $\varPhi_{K_S}$ of $S$ contracts no 
$(-2)$-curve except those contracted by $S \to \mathcal{C}$. 
Thus $\mathcal{C}$ has no $(-1)$-curve contained in the image 
of a fundamental cycle of $S$.  \qed

Note that the surfaces $S$'s in the Proposition above are  
topologically simply connected.  
For regular surfaces with $c_1^2 = 8$ and $p_g = 4$ with 
non trivial torsion, see \cite{on2chi-2}. 
In a sequel to this paper, 
the example in Proposition \ref{prop:bauerpignatellim0} 
will be further studied, being deformed to surfaces with 
non-hyperelliptic genus $3$ fibrations.


\begin{thebibliography}{11}

\bibitem{akgloballocal}
{\sc Ashikaga,~T., Konno,~K.} 
Global and local properties of pencils of algebraic curves, Algebraic geometry 2000, Azumino (Hotaka), 
{\em Adv. Stud. Pure Math.,36, Math. Soc. Japan, Tokyo}, {\bf } (2002), 1--49.

\bibitem{complexsurf}
{\sc Barth,~W., Peters,~C., {Van de Ven},~A.} 
Compact Complex Surfaces, 
{\em Springer-Verlag}, (1984). 

\bibitem{pg4c7bauer}
{\sc Bauer,~I.} 
Surfaces with $K^2=7$ and $p_g = 4$, 
{\em Mem.\, Amer.\, Math.\, Soc.}, {\bf 152} (2001), no.\, 721.

\bibitem{caninvpg4c8}
{\sc Bauer,~I., Pignatelli,~R.} 
Surfaces with $K^2=8$, $p_g = 4$ and canonical involution, 
{\em Osaka J. Math.}, {\bf 46} (2009), no.\, 3, 799--820.

\bibitem{catonmodulijdg}
{\sc Catanese,~F.} 
On the moduli spaces of surfaces of general type, 
{\em J.\, Differential Geom.}, {\bf 19} (1984), no.\, 2, 483--515.

\bibitem{pgq1catcil}
{\sc Catanese,~F., Ciliberto,~C.} 
Surfaces with $p_g = q = 1$. Problems in the theory of surfaces and their classification (Cortona, 1988), 
{\em Sympos. Math., XXXII, Academic Press, London}, {\bf } (1991), 49--79.

\bibitem{catcilsymm}
{\sc Catanese,~F., Ciliberto,~C.} 
Symmetric products of elliptic curves and surfaces of general type with $p_g = q = 1$, 
{\em J. Algebraic Geom.}, {\bf 2} (1993), 389--411.


\bibitem{catlipign}
{\sc Catanese,~F., Liu,~W., Pignatelli,~R.} 
The moduli space of even surfaces of general type with $K^2 = 8$, $p_g = 4$ and $q = 0$, 
{\tt arXiv:1209.0034 [math.AG]}, {\bf } (2012),.


\bibitem{fibrationsI'}
{\sc Catanese,~F., Pignatelli,~R.} 
Fibrations of low genus, I, 
{\em Ann. Sci. \'{E}cole Norm. Sup. (4)}, {\bf 39} (2006), 1011--1049.

\bibitem{cilcansfpg4}
{\sc Ciliberto,~C.} 
Canonical surfaces with $p_g=p_a=4$ and $K^2=5, \ldots , 10$, 
{\em Duke Math. J.}, {\bf 48} (1981), no.\, 1, 121--157.

\bibitem{on2chi-2}
{\sc Ciliberto,~C., Mendes Lopes,~M.} 
On surfaces of general type with $K^2 = 2\chi -2$ and non-trivial torsion, 
{\em Geom. Dedicata}, {\bf 66} (1997), no. 3, 313--329.

\bibitem{enriquessuperfici}
{\sc Enriques,~F.} 
Le superficie algebriche, 
{\em Nicola Zanichelli, Bologna}, {\bf } (1949).

\bibitem{fujitakaehler}
{\sc Fujita,~T.} 
On K\"{a}hler fiber spaces over curves, 
{\em J. Math. Soc. Japan}, {\bf 30} (1978), 779--794.

\bibitem{smallc1-1}
{\sc Horikawa,~E.} 
Algebraic surfaces of general type with small $c_1^2$. I, 
{\em Ann. of Math.}, {\bf 104} (1976), no.\ 2, 358--387.

\bibitem{smallc1-2}
{\sc Horikawa,~E.} 
Algebraic surfaces of general type with small $c_1^2$. II, 
{\em Invent. Math.}, {\bf 37} (1976), 121--155.

\bibitem{horikawapencils}
{\sc Horikawa,~E.} 
On algebraic surfaces with pencils of curves of genus 2, 
{\em Complex analysis and algebraic geometry, Iwanami Shoten, Tokyo}, {\bf } (1977), 79--90.

\bibitem{kodoncompact2}
{\sc Kodaira,~K.} 
On compact analytic surfaces. II, 
{\em Ann. of Math. (2)}, {\bf 77} (1963), 563--626.

\bibitem{kodoncompact3}
{\sc Kodaira,~K.} 
On compact analytic surfaces. III, 
{\em Ann. of Math. (2)}, {\bf 78} (1963), 1--40

\bibitem{cliffkonno}
{\sc Konno,~K.} 
Clifford index and the slope of fibered surfaces, 
{\em J. Algebraic Geom.}, {\bf 8} (1999), 207–-220.

\bibitem{1-2-3}
{\sc Konno,~K.} 
1-2-3 theorem for curves on algebraic surfaces, 
{\em J. Reine. Angew. Math.}, {\bf 533} (2001), 171--205.

\bibitem{questionofreid}
{\sc Konno,~K., Mendes Lopes,~M.} 
On a question of Miles Reid, 
{\em Manuscripta Math.}, {\bf 100} (1999), no.\ 1, 81--86.

\bibitem{thesismlopes}
{\sc Mendes Lopes,~M.} 
The relative canonical algebra for genus $3$ fibrations, 
{\em Thesis, Warwick University}, {\bf } (1988).

\bibitem{mlopespardini3chi}
{\sc Mendes Lopes,~M.，Pardini,~R.} 
On the algebraic fundamental group of surfaces with $K^2 \leq 3 \chi$, 
{\em J.\, Differential Geom.}, {\bf 77} (2007), no.\, 2, 189--199.

\bibitem{namuengenus2}
{\sc Namikawa,~Y., Ueno,~K.} 
The complete classification of fibres in pencils of curves of genus two, 
{\em Manuscripta Math.}, {\bf 9} (1973), 143--186.

\bibitem{ogggenus2}
{\sc Ogg,~A.~P.} 
On pencils of curves of genus two, 
{\em Topology}, {\bf 5} (1966), 355--362.


\bibitem{oliverioevenpg4c8}
{\sc Oliverio,~P.} 
On even surfaces of general type with $K^2=8$, $p_g=4$, $q=0$, 
{\em Rend.\, Sem.\, Mat.\, Univ.\, Padova}, {\bf 113} (2005), 1--14.

\bibitem{reidpencils}
{\sc Reid,~M.} 
Problems on pencils of small genus, 
{\em preprint}, {\bf } (1990).

\bibitem{supinopg4c8}
{\sc Supino,~P.} 
On moduli of regular surfaces with $K^2 = 8$ and $p_g = 4$, 
{\em Port.\, Math.\, (N.S.)}, {\bf 60} (2003) no.\, 3, 353--358.

\bibitem{xiaogenre2}
{\sc Xiao,~G.} 
Surfaces fibr\'{e}es en courbes de genre deux. Lecture Notes in Mathematics, 1137, 
{\em Springer-Verlag, Berlin}, {\bf } (1985).


\bibitem{pi1hyperell}
{\sc Xiao,~G.} 
$\pi_1$ of elliptic and hyperelliptic surfaces, 
{\em Internat. J. math.}, {\bf 2} (1991) 599--615.


\end{thebibliography}

{\sc 
Masaaki Murakami, 

University of Bayreuth, Lehrstuhl Mathematik VIII, 
Universitaetsstrasse 30, 
D-95447 Bayreuth, Germany}

{\it E-mail address}:\, \texttt{Masaaki.Murakami@uni-bayreuth.de}


\end{document}